\documentclass[
11pt, a4paper, USenglish
]{article}

\usepackage{import}

\usepackage[USenglish]{babel}
\usepackage{import}

\usepackage[utf8]{inputenc}

\usepackage[T1]{fontenc}
\usepackage{cite}

\usepackage{mathscinet}

\usepackage{amsmath, amscd, amsthm}
\usepackage{thmtools}
\usepackage{amssymb}
\usepackage{mathtools}

\usepackage{tikz-cd}
\usepackage{tikz}

\usepackage{ifthen}

\newcommand{\ifisemptythenelse}[3]{
  \if\relax\detokenize{#1}\relax
    #2%
  \else
    #3%
  \fi
}

\usepackage{dsfont}
\usepackage[mathscr]{eucal}

\newcommand\ol{\overline}

\newcommand{\mathblackboardboldnew}{\mathds}

\newcommand{\BN}{\mathblackboardboldnew N}

\newcommand{\BZ}{\mathblackboardboldnew Z}

\newcommand{\newmathcal}{\mathcal}

\newcommand{\calB}{\newmathcal{B}}
\newcommand{\calC}{\newmathcal{C}}
\newcommand{\calD}{\newmathcal{D}}

\newcommand{\calF}{\newmathcal{F}}

\newcommand{\calJ}{\newmathcal{J}}

\newcommand{\calP}{\newmathcal{P}}

\newcommand{\calT}{\newmathcal{T}}

\newcommand{\calX}{\newmathcal{X}}
\newcommand{\calY}{\newmathcal{Y}}

\newcommand{\rmK}{\mathrm{K}}

\newcommand{\rmS}{\mathrm{S}}

\newcommand{\qquote}[1]{``#1''}
\newcommand{\roughly}[1]{``#1''}

\newcommand{\introduce}{\textbf}
\newcommand{\latin}{\emph}
\newcommand{\buzzword}{\emph}

\newcommand{\pureword}[2]{
  \ifisemptythenelse%
  {#2}%
  {#1}%
  {\ErrorArgumentNotEmpty}%
}

\newcommand{\Ktheory}[1]{\pureword{$\rmK$-theory}{#1}}

\newcommand{\Sconstruction}[1]{\pureword{$\rmS$-construction}{#1}}

\newcommand{\inftycat}[1]{\pureword{$\infty$-category}{#1}}
\newcommand{\inftycats}[1]{\pureword{$\infty$-categories}{#1}}

\newcommand{\inftyoperads}[1]{\pureword{$\infty$-operads}{#1}}

\newcommand{\resp}{{resp.\ }}
\newcommand\auxiliarytoggleablelatin\latin

\newcommand{\ie}{\auxiliarytoggleablelatin{i.e.,} }

\newcommand{\twonames}[2]{{#1} and {#2}}
\newcommand{\threenames}[3]{{#1}, {#2} and {#3}}
\newcommand{\fivenames}[5]{{#1}, {#2}, {#3}, {#4} and {#5}}

\newcommand{\Cech}[1]{\pureword{\v{C}ech}{#1}}
\newcommand{\nameDK}[1]{\pureword{\twonames{Dyckerhoff}{Kapranov}}{#1}}
\newcommand{\Galvez}[1]{\pureword{G{\'a}lvez}{#1}}
\newcommand{\nameGCKT}[1]{\pureword{\threenames{\Galvez{}-Carrillo}{Kock}{Tonks}}{#1}}
\newcommand{\nameBOORS}[1]{\pureword{\fivenames{Bergner}{Osorno}{Ozornova}{Rovelli}{Scheimbauer}}{#1}}

\newcommand{\ConnesCC}[1]{\pureword{Connes' cyclic category}{#1}}

\newcommand{\TDyckerhoff}[1]{\pureword{T.~Dyckerhoff}{#1}}
\newcommand{\GJasso}[1]{\pureword{G.~Jasso}{#1}}

\newcommand{\theauthor}[1]{\pureword{the author}{#1}} 

\newtheoremstyle{newtheorem}{}{}{}{}{\bfseries}

\DeclareSymbolFont{extraup}{U}{zavm}{m}{n}
\DeclareMathSymbol{\varheart}{\mathalpha}{extraup}{86}
\DeclareMathSymbol{\vardiamond}{\mathalpha}{extraup}{87}
\DeclareMathSymbol{\varclub}{\mathalpha}{extraup}{84}
\DeclareMathSymbol{\varspade}{\mathalpha}{extraup}{85}

\newcommand{\somberend}{\ensuremath{\Diamond}}

\ifx\fancyends\undefined
\newcommand{\theoremend}{\ensuremath{\square}}
\newcommand{\definitionend}{\somberend}
\newcommand{\exampleend}{\somberend}
\newcommand{\questionend}{\somberend}
\newcommand{\proofend}{\ensuremath{\blacksquare}}
\else
\newcommand{\theoremend}{\ensuremath{\varheart}}
\newcommand{\definitionend}{\ensuremath{\clubsuit}}
\newcommand{\exampleend}{\ensuremath{\diamondsuit}}
\newcommand{\questionend}{\ensuremath{\spadesuit}}
\newcommand{\proofend}{\ensuremath{\square}}
\fi

\theoremstyle{definition}

\newtheorem{Masterthm}{Masterthm}[subsection]

\declaretheorem[name=Construction ,style=definition,qed={\definitionend},sibling=Masterthm]{Cstr}
\declaretheorem[name=Construction ,style=definition,qed={\definitionend},unnumbered]{Cstr*}

\declaretheorem[name=Corollary,style=definition,qed={\theoremend},sibling=Masterthm]{Cor}
\declaretheorem[name=Definition,style=definition,qed={\definitionend},sibling=Masterthm]{Def}
\declaretheorem[name=Definition,style=definition,qed={\definitionend},unnumbered]{Def*}

\declaretheorem[name=Goal,style=definition,qed={\questionend},unnumbered]{Goal*}

\declaretheorem[name=Lemma,style=definition,qed={\theoremend},sibling=Masterthm]{Lem}

\declaretheorem[name=Observation,style=definition,qed={\theoremend},sibling=Masterthm]{Obs}
\declaretheorem[name=Proof,style=definition,qed={\proofend}, numbered=no]{Prf}

\declaretheorem[name=Proposition,style=definition,qed={\theoremend},sibling=Masterthm]{Prop}
\declaretheorem[name=Question,style=definition,qed={\questionend},sibling=Masterthm]{Qstn}

\declaretheorem[name=Theorem,style=definition,qed={\theoremend},sibling=Masterthm]{Thm}
\declaretheorem[name=Theorem,style=definition,qed={\theoremend},unnumbered]{Thm*}

\theoremstyle{remark}
\declaretheorem[name=Remark ,style=remark,qed={\exampleend},sibling=Masterthm]{Rem}

\declaretheorem[name=Fact,style=remark,qed={\exampleend},unnumbered]{Fact*}

\declaretheorem[name=Example ,style=remark,qed={\exampleend},sibling=Masterthm]{Expl}


\numberwithin{equation}{section}

\declaretheorem[name=Theorem, style=definition, qed={\somberend}]{Theorem}
\declaretheorem[name=Theorem, style=definition, qed={\somberend},unnumbered]{Theorem*}
\declaretheorem[name=Corollary, style=definition, qed={\somberend}]{Corollary}
\declaretheorem[name=Corollary, style=definition, qed={\somberend},unnumbered]{Corollary*}

\declaretheorem[name=Observation, style=definition, qed={\somberend}, unnumbered]{Observation*}

\declaretheorem[name=Conjecture, style=definition, qed={\somberend},unnumbered]{Conjecture*}

\declaretheorem[name=Definition, style=definition, qed={\somberend}, unnumbered]{Definition*}

\declaretheorem[name=Remark ,style=remark,qed={\somberend},unnumbered]{Remark*}

\ifx\nopagestyle\undefined
\usepackage{geometry}
\usepackage{enumitem}

\setlist{itemsep=-3pt, topsep=1pt}
\setlist[enumerate, 1]{label=(\arabic*), ref=(\arabic*)}

\newcommand{\axiomlabelstyle}[1]{(\bfseries{#1}\arabic*)} 

\usepackage{footnote}

\def\defaultdepthtoc{3}

\usepackage[titletoc]{appendix}				%
\usepackage[nottoc,notlot,notlof]{tocbibind} 	
\addto\captionsUSenglish{
  \renewcommand{\contentsname}%
    {Contents:}%
}
\addtocontents{toc}{\protect\setcounter{tocdepth}{\defaultdepthtoc}}

\usepackage{abstract} 

\usepackage{setspace}

\usepackage{fancyhdr}
\usepackage{lastpage}

\fi

\newlength{\longarrow}
\newlength{\arrow}
\newlength{\blablaone}
\newlength{\blablatwo}

\newcommand{\lra}{\longrightarrow}
\newcommand{\lla}{\longleftarrow}

\newcommand{\hra}{\hookrightarrow}
\newcommand{\lhra}{\xhra{}}

\newcommand{\xhra}[2][]{
	\settowidth{\blablaone}{\scriptsize$#1$}
	\settowidth{\blablatwo}{\scriptsize$#2$}
	\pgfmathsetlength{\arrow}{max(\blablaone,\blablatwo,\longarrow)}
	\xhookrightarrow[{\mathmakebox[\arrow]{#1}}]{\mathmakebox[\arrow]{#2}}%
}
\newcommand{\xra}[2][]{
	\settowidth{\blablaone}{\scriptsize$#1$}
	\settowidth{\blablatwo}{\scriptsize$#2$}
	\pgfmathsetlength{\arrow}{max(\blablaone,\blablatwo,\longarrow)}
	\xrightarrow[{\mathmakebox[\arrow]{#1}}]{\mathmakebox[\arrow]{#2}}%
}

\newcommand{\xlra}[2][]{
	\settowidth{\blablaone}{\scriptsize$#1$}
	\settowidth{\blablatwo}{\scriptsize$#2$}
	\pgfmathsetlength{\arrow}{max(\blablaone,\blablatwo,\longarrow)}
	\xleftrightarrow[{\mathmakebox[\arrow]{#1}}]{\mathmakebox[\arrow]{#2}}%
}

\newcommand{\lmapsto}{
	\xmapsto{\mathmakebox[\longarrow]{}}%
}

\newcommand{\lrlas}{
  \mathrel{\substack{\lra \\[-.65ex] \lla}}
}

\newcommand\adjarrows{\lrlas}
\newcommand\ladjarrows\adjarrows
\newcommand\radjarrows\adjarrows

\usepackage{xstring}

\newcommand\intxt[1]{\qquad\text{#1}\qquad}

\newcommand\isCartesian{\drar[phantom, description, very near start,"\lrcorner"]}
\newcommand\iscoCartesian{\drar[phantom, description, very near end,"\ulcorner"]}
\newcommand\isbiCartesian{\drar[phantom, description, "\square"]}

\newcommand\cdsquare[9][]{%
  \def\obja{#2}%
  \def\objb{#3}%
  \def\objc{#4}%
  \def\objd{#5}%
  \def\mora{#6}
  \def\morb{#7}
  \def\morc{#8}
  \def\mord{#9}
  \begin{tikzcd}[ampersand replacement=\&]
    \obja\ar[r,"{\mora}"]\ar[d,"{\morb}"']%
    \IfEqCase{#1}{
      {C}{\isCartesian}
      {cC}{\iscoCartesian}
      {bC}{\isbiCartesian}
    }
    \&%
    \objb\ar[d,"{\morc}"]%
    \\%
    \objc\ar[r,"{\mord}"']%
    \&%
    \objd%
  \end{tikzcd}%
}

\newcommand\cdcubeNA[9][]{%
  \def\obja{#2}%
  \def\objb{#3}%
  \def\objc{#4}%
  \def\objd{#5}%
  \def\obje{#6}
  \def\objf{#7}
  \def\objg{#8}
  \def\objh{#9}
  \begin{tikzcd}[ampersand replacement=\&]
    \obja\ar[rr]\ar[dr]\ar[dd]\&
    \&\objb\ar[dr]\ar[dd]\&\\
    \&\objc\ar[rr,crossing over]\&
    \&\objd\ar[dd]\\
    \obje\ar[rr]\ar[dr]\&
    \&\objf\ar[dr]\&\\
    \&\objg\ar[rr]\ar[from=uu,crossing over]\&
    \&\objh\\
  \end{tikzcd}%
}

\newcommand\cdsquareNA[5][]{\cdsquare[#1]{#2}{#3}{#4}{#5}{}{}{}{}}

\DeclareMathOperator*{\colim}{colim}

\DeclareMathOperator{\Aut}{Aut}

\DeclareMathOperator{\Fun}{Fun}

\newcommand\Id{\mathrm{Id}}

\newcommand\blank{-}

\newcommand\inv{{-1}}
\newcommand\op{\mathrm{op}}

\newcommand\basepoint\star
\newcommand\bigdisjunion{\dot\bigcup}
\newcommand\card[1]{\left|#1\right|}
\newcommand\disjunion{\mathbin{\dot\cup}}
\newcommand\emptytuple\varnothing
\newcommand\invimage[2]{{#1}^{\inv}#2}
\newcommand\powerset[1]{\calP(#1)}
\newcommand\powersetgeop[2]{\calP^\op_{\geq #2}(#1)}
\newcommand\powersetleop[2]{\calP^\op_{\leq #2}(#1)}
\newcommand\powersetne[1]{\calP_{*}(#1)}
\newcommand\powersetneop[1]{\calP^\op_{*}(#1)}
\newcommand\powersetop[1]{\calP^\op(#1)}
\newcommand\preimage[2]{\invimage{#1}{\set{#2}}}
\newcommand\setmid{\,\middle|\,}
\newcommand\setP[2]{\left\{#1\setmid #2\right\}}
\newcommand\set[1]{\left\{#1\right\}}
\newcommand\tupleP[2]{\left(#1\setmid#2\right)}
\renewcommand\emptyset\varnothing

\newcommand\cygrp[1]{\rquot \BZ {#1}}

\newcommand\const[1][]{\ifisemptythenelse{#1}{\mathrm{const}}{\mathrm{const}_{#1}}}

\newcommand\fiberproduct[1]{\times_{#1}}
\newcommand\ladjto\dashv
\newcommand\radjto\vdash
\newcommand\overcat[2]{{{#1}_{/#2}}}

 \newcommand\rquot[2]{
        \mathchoice
            {
                \text{\raise1ex\hbox{${#1}$}}\Big/\lower1ex\hbox{${#2}$}%
            }
            {
                #1\,/\,#2
            }
            {
                #1/#2
            }
            {
                #1/#2
            }
    }
\newcommand{\doubleslash}{\mathbin{
  \mathchoice{\Big/\mkern-10mu\Big/}
    {/\mkern-6mu/}
    {/\mkern-5mu/}
    {/\mkern-5mu/}}}

\newcommand{\qquot}[2]{
\mathchoice
	    {
                \text{\raise1ex\hbox{${#1}$}}\doubleslash\lower1ex\hbox{${#2}$}%
            }
	{{#1}\doubleslash{#2}}
	{{#1}\doubleslash{#2}}
	{{#1}\doubleslash{#2}}
} 
    
\newcommand\lquot[2]{
        \mathchoice
            {
                \lower1ex\hbox{${#2}$}\Big\backslash\text{\raise1ex\hbox{${#1}$}}%
            }
            {
                #2\setminus#1
            }
            {
                #2\setminus#1
            }
            {
                #1\setminus#2
            }
    }
    
 \newcommand\lrquot[3]{
        \mathchoice
            {
                \lower1ex\hbox{${#2}$}\Big\backslash\text{\raise1ex\hbox{${#1}$}}\Big/\lower1ex\hbox{${#3}$}%
            }
            {
                #2\backslash#1/#3
            }
            {
                #2\backslash#1/#2
            }
            {
                #1\backslash#2/#2
            }
    }

\newcommand\restr[3][]{{
  \left.\kern-\nulldelimiterspace 
  {#2} 
  \vphantom{\big|} 
  \right|_{#3}^{#1} 
  }}

\newcommand{\category}[1]{\mathrm{\mathbf{#1}}}
\newcommand{\cat}\category

\newcommand{\Set}{\cat{Set}}
\newcommand{\poSet}{\cat{Pos}}

\newcommand\0{{\numD 0}}
\newcommand{\m}{{\numD m}}
\newcommand{\n}{{\numD n}}
\newcommand\numD[1]{{[#1]}}
\newcommand\numL[1]{{\langle#1\rangle}}

\newcommand\Deltaug{\Delta_+}

\newcommand\Deltaleq[1]{\Delta_{\leq #1}}
\newcommand\Deltall[1]{\Delta_{< #1}}
\newcommand{\Dop}{{\Delta^{\op}}}

\newcommand{\cofacemap}[1]{\mathrm{d}^{#1}}

\newcommand{\codegmap}[1]{\mathrm{s}^{#1}}

\newcommand\join{\mathbin\star}
\DeclareMathOperator*{\bigjoin}{\bigstar}
\newcommand\bigjoinl{\bigjoin\limits}
\newcommand\lpath{{P^\triangleleft}}
\newcommand\rpath{{P^\triangleright}}

\newcommand\Segalspine[1]{\mathrm{Sp}[#1]}
\newcommand{\simplex}[1]{\Delta^{#1}}

\usepackage{hyperref}

\addto\extrasUSenglish{%
}

\usepackage[shortcuts]{extdash}

\begin{document}

\newcommand{\nameBdBW}[1]{\pureword{\twonames{Boavida de Brito}{Weiss}}{#1}}

\title{Higher Segal spaces via higher excision}
\def\myauthor{Tashi Walde}
\author{\myauthor\thanks{Rheinische Friedrich-Wilhelms-Universität Bonn,
    \url{http://www.math.uni-bonn.de/people/twalde/}
}}

\date{}
\maketitle

\begin{abstract}
  We show that the various higher Segal conditions of \nameDK{}
  can all be characterized in purely categorical terms
  by higher excision conditions
  (in the spirit of Goodwillie--Weiss manifold calculus)
  on the simplex category $\Delta$ and the cyclic category $\Lambda$.
\end{abstract}

\tableofcontents

\newpage
\newcommand{\localC}{\calC}
\newcommand{\localD}{\calD}
\newcommand{\localX}{\calX}
\newcommand{\xdc}{\localX\colon\Dop\to\localC}
\newcommand{\excube}{Q}
\newcommand{\ff}{f}
\newcommand{\II}{I}

\newcommand{\clawon}{\models}
\newcommand\exclaw{\calF}
\newcommand\refines{\preceq}
\newcommand\Cechcube[1]{\check{\mathrm C}#1}

\newcommand{\Xonsset}[1]{\localX_{#1}}
\newcommand{\exsset}{K}

\newcommand{\canactivemap}[1]{a_{#1}}
\newcommand{\Dactop}{\Delta^{\mathrm{act},\op}}
\newcommand{\Dact}{\Delta^{\mathrm{act}}}
\newcommand{\Dlactop}{\Delta^{\mathrm{lact},\op}}
\newcommand{\Dlact}{\Delta^{\mathrm{lact}}}
\newcommand{\Dractop}{\Delta^{\mathrm{ract},\op}}
\newcommand{\Dract}{\Delta^{\mathrm{ract}}}
\newcommand{\Drstr}{\Delta^{\mathrm{rstr}}}
\newcommand{\Deltainfty}{\Delta_\infty}
\newcommand{\Dactgeq}[1]{\Dact_{\geq{#1}}}
\newcommand{\Dactgeqop}[1]{\Dactop_{\geq{#1}}}
\newcommand{\numDminus}[2]{\numD{#1\setminus #2}}
\newcommand{\readbl}[1]{{{#1}^+}}
\newcommand{\readbr}[1]{{{#1}^-}}

\newcommand{\ssetofprecover}[1]{\widetilde{#1}}
\newcommand{\reducedclaw}[1]{\ssetofprecover{#1}}

\newcommand{\kk}{k}
\newcommand{\dd}{d}
\newcommand\lSegalclaw[2]{\operatorname{lSeg}^{#1}_{#2}}

\newcommand{\evenSegal}[1]{$2{#1}$-Segal}
\newcommand{\loweroddSegal}[1]{lower ($2{#1}-1$)-Segal}

\section{Introduction}

In his seminal 2001 paper~\cite{Rezk2001}, Rezk introduced
the notion of Segal objects in order to
describe monoids (or, more generally, categories)
which are not strictly associative
but only associative up to a coherent system of higher homotopies.
A Segal object in an ($\infty$-)category $\localC$
is a simplicial object in $\localC$---%
\ie a functor $\localX\colon\Dop\to\localC$ on the simplex category $\Delta$---%
satisfying a certain family of descent conditions.

The starting point for this work is the
easy but little-known observation
that Segal objects can be characterized
by a condition which is purely categorical,
in the sense that it can be defined without having to know anything
about the inner workings of $\Delta$.

\begin{Observation*}
  A simplicial object $\Dop\to\localC$ is Segal if and only if
  it sends biCartesian squares in $\Delta$ to Cartesian squares in $\localC$.
\end{Observation*}

\newcommand\cycPol[2]{\mathrm{C}(#2,#1)} 
\newcommand{\extriang}{\calT}

In 2012, \nameDK{} generalized Rezk's Segal condition
and introduced what they call \buzzword{higher Segal spaces}.
Their definition is very geometric in nature:
They consider the so called
\buzzword{cyclic polytopes} $\cycPol{\dd}{n}$,
defined as the convex hull
of $n+1$ points on the $\dd$-dimensional moment curve
$t\mapsto(t,t^2,\dots,t^\dd)$.
The main feature of these polytopes in this context
is that they have two canonical triangulations,
called the \buzzword{lower triangulation}
and the \buzzword{upper triangulation}, respectively.
Each of these triangulations defines a simplicial subcomplex $\extriang$
of the standard $n$-simplex $\simplex n$;
\nameDK{} then impose conditions on simplicial objects
by requiring that the value\footnote{
  Every simplicial object
 can be canonically evaluated on simplicial sets
 by Kan extension along the Yoneda embedding;
 see \autoref{sec:membrane-spaces}.
} on
the inclusion $\extriang\hra\simplex n$ is an equivalence:
a simplicial object is called
\introduce{lower} (\resp \introduce{upper}) \introduce{$\dd$-Segal}
if this is true for the lower (\resp upper) triangulation of $\cycPol{\dd}{n}$
and \introduce{$\dd$-Segal}
if this is true for all triangulations of $\cycPol{\dd}{n}$.\\

The purpose of this article is to characterize
the various flavors of higher Segal conditions
in terms of purely categorical notions of higher excision.
We first do this for \loweroddSegal{\kk} spaces,
since they are in a precise sense\footnote{
  This vague assertion is made precise by the
  path space criterion~\cite[Proposition~2.7]{Poguntke2017}
  which expresses all higher Segal conditions in terms
  of \loweroddSegal{\kk} conditions.
}
the most fundamental amongst all versions of higher Segal spaces.
The following is the first main result of this paper:

\begin{Theorem}[\autoref{thm:main}]
  \label{thm:intro-main-odd}
  Let $\xdc$ be a simplicial object in an $\infty$-category $\localC$
  with finite limits.
  The following are equivalent:
  \begin{enumerate}
  \item
    the simplicial object $\localX$ is \loweroddSegal{\kk};
  \item
    \label{item:12}
    the functor $\localX$ sends every strongly biCartesian\footnote{
      A cube is strongly biCartesian
      if each of its $2$-dimensional faces is biCartesian;
      see \autoref{def:strongly-biCart}.
    } ($k+1$)-dimensional cube in $\Delta$ to a limit diagram in $\localC$.
    \qedhere
  \end{enumerate}
\end{Theorem}

We call a functor $\localD^\op\to\localC$ satisfying condition~\ref{item:12}
of \autoref{thm:intro-main-odd} \introduce{weakly $\kk$-excisive};
compare this with Goodwillie's calculus of functors~\cite{Goodwillie91},
where a (covariant) functor $\localD\to\localC$ is called $\kk$-excisive
if it sends strongly \emph{co}Cartesian ($\kk+1$)-dimensional cubes in $\localD$
to limit diagrams in $\localC$.\\

We illustrate \autoref{thm:intro-main-odd} with some examples.

\begin{figure}[t]
  \centering
  \includegraphics{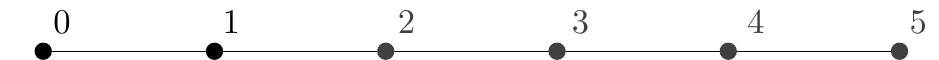}
  \caption[The lower $1$-Segal triangulation]{
    The lower triangulation of the
    cyclic polytope $\cycPol{1}{n}$, here depicted with $n=5$.
  }
  \label{fig:cyc-1-n}
\end{figure}

\begin{figure}[t]
  \centering
  \includegraphics{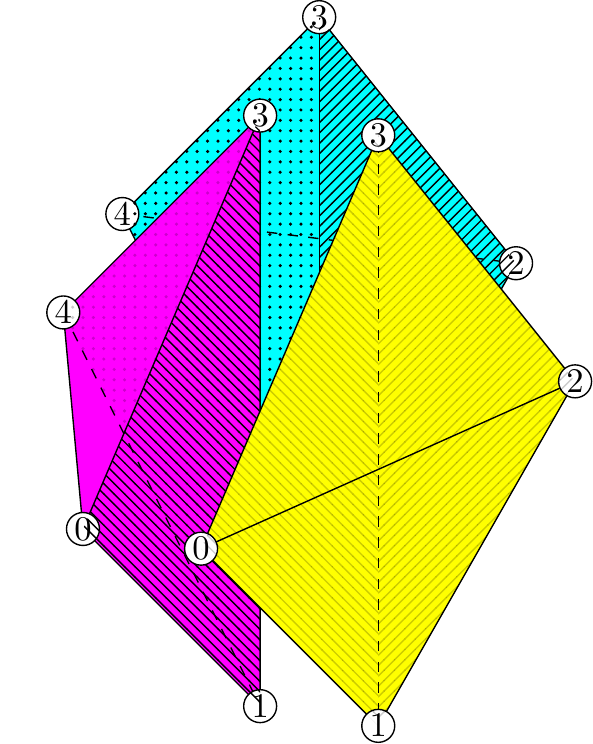}
  \caption[The lower triangulation of $\cycPol{3}{4}$]{
    The three $3$-simplices
    $\simplex{\set{1234}}$, $\simplex{\set{0134}}$ and $\simplex{\set{0123}}$
    (depicted in cyan, magenta and yellow, respectively)
    assemble into the lower triangulation
    of the double triangular pyramid $\cycPol{3}{4}$.
  }
  \label{fig:cyc-3-4}
\end{figure}

\begin{itemize}
\item
  The cyclic polytope $\cycPol{1}{n}$ is just the interval
  $\simplex{\set{0,n}}$;
  its lower triangulation (see \autoref{fig:cyc-1-n})
  yields the simplicial complex
  \begin{equation}
    \Segalspine{n}\coloneqq
    \simplex{\set{0,1}}
    \cup
    \dots
    \cup
    \simplex{\set{n-1,n}}
    \subset\simplex n.
  \end{equation}
  Rezk's Segal condition for a simplicial object says precisely that the
  inclusion $\Segalspine{n}\hra\simplex{n}$ needs to be sent to an
  equivalence; this is what \nameDK{} call the lower $1$-Segal condition.
  For $n=1$, this condition says precisely that the biCartesian square
  \begin{equation}
    \label{eq:intro-min-1Segal-square}
    \cdsquareNA[bC]{1}{12}{01}{012}
  \end{equation}
  in $\Delta$ needs to be sent to a limit diagram.
  More generally, every square of the form
  \begin{equation}
    \cdsquareNA[bC]
    {\set{i}}{\set{i,\dots,n}}
    {\set{0,\dots,i}}{\set{0,\dots,n}}
  \end{equation}
  (for $0<i<n$) is biCartesian in $\Delta$;
  it is in fact an often used characterization of Segal objects
  to require these squares to be sent to pullbacks.
\item
  The cyclic polytope $\cycPol{3}{4}$ is a double triangular pyramid;
  its lower triangulation
  (see \autoref{fig:cyc-3-4})
  induces the simplicial complex
  \begin{equation}
  \extriang=
  \simplex{\set{1,2,3,4}}
  \cup
  \simplex{\set{0,1,3,4}}
  \cup
  \simplex{\set{0,1,2,3}}\subset\simplex 4.
\end{equation}
  By definition, a simplicial object satisfies
  the first lower $3$-Segal condition
  if it sends the canonical inclusion
  $\extriang\hra\simplex 4$ to an equivalence;
  this is equivalent to sending the cube
  \begin{equation}
    \label{eq:6}
    \cdcubeNA
    {13}{134}
    {123}{1234}
    {013}{0134}
    {0123}{01234}
  \end{equation}
  which is strongly biCartesian in $\Delta$, to a limit diagram.
\end{itemize}

In general, the first non-trivial \loweroddSegal{\kk} condition
(\ie the one for $n=2\kk$)
can always be expressed in terms of a strongly biCartesian cube in $\Delta$
of dimension $k+1$ and this cube is the unique such cube which is
in a certain sense \roughly{basic}.
However, for bigger $n$
both the number of simplices in the lower triangulation of $\cycPol{2\kk-1}{n}$
and the number of basic strongly biCartesian cubes grows very rapidly
so that, \latin{a priori},
the behavior of weakly $\kk$-excisive simplicial objects
and \loweroddSegal{\kk} objects diverges dramatically.\\

Since the introduction of higher Segal spaces,
most interest in the area was garnered by $2$-Segal spaces;
more precisely by $2$-Segal spaces that satisfy
an additional condition called \buzzword{unitality}.
For example, unital $2$-Segal spaces were studied
by Dyckerhoff from the perspective of Hall algebras~\cite{Dyckerhoff2018}
and by \nameGCKT{}~\cite{GCKT2018a,GCKT2018b,GCKT2018c}
from the perspective of bialgebras arising in
combinatorics\footnote{
  Unital $2$-Segal spaces are called
  \buzzword{decomposition spaces} in this context.
}.
The \inftycat{} of unital $2$-Segal spaces
was identified by \theauthor{}~\cite{Walde2017}
as a certain sub-\inftycat{} of the \inftycat{} of \inftyoperads{}
and recently by Stern~\cite{Stern2019}
as a certain \inftycat{} of algebras in correspondences.
The main source of examples for unital $2$-Segal objects is
Waldhausen's \Sconstruction{} from algebraic \Ktheory{};
\nameBOORS{}~\cite{BOORS2018}
showed that in a certain sense every unital $2$-Segal
space arises this way;
Poguntke~\cite{Poguntke2017}
generalized Waldhausen's construction to higher dimensions,
thus providing many examples for \evenSegal{\kk} spaces.
Furthermore, \emph{cyclic} unital $2$-Segal spaces---which
can be identified with certain cyclic \inftyoperads{}~\cite{Walde2017}
or with Calabi-Yau algebras in correspondences~\cite{Stern2019}---
play a central role in the
construction of topological Fukaya categories of marked surfaces
due to \nameDK{}~\cite{DyckerhoffKapranov2018}.

We show that $2$-Segal spaces,
and more generally \evenSegal{\kk} spaces,
can be characterized by a relative version of higher weak excision
which involves \ConnesCC{} $\Lambda$.

\begin{Theorem}[\autoref{thm:main}]
  \label{thm:intro-main-even}
  Let $\xdc$ be a simplicial object in an $\infty$-category $\localC$
  with finite limits.
  The following are equivalent:
  \begin{enumerate}
  \item
    the simplicial object $\localX$ is \evenSegal{\kk};
  \item
    the functor $\localX$ sends to Cartesian cubes in $\localC$
    those ($\kk+1$)-dimensional cubes in $\Delta$
    which become strongly biCartesian in $\Lambda$
    (under the canonical functor $\Delta\to\Lambda$).
    \qedhere
  \end{enumerate}
\end{Theorem}

\begin{figure}[t]
  \centering
  \includegraphics{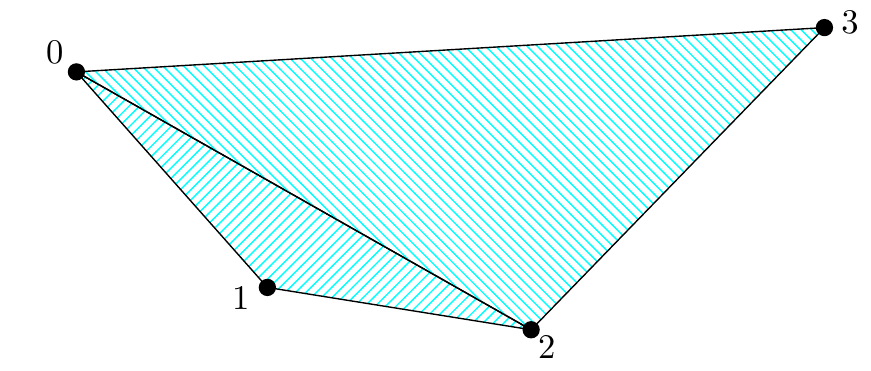}
  \includegraphics{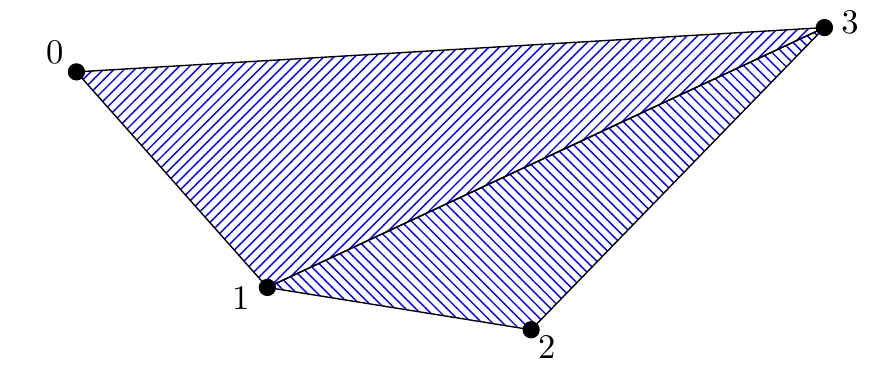}
  \caption[Triangulations of $\cycPol{2}{4}$]{
    The lower and the upper
    triangulations of the cyclic polytope $\cycPol{2}{4}$.
  }
  \label{fig:cyc-2-4} 
\end{figure}

We again illustrate the theorem with some examples:

\begin{itemize}
\item
  The square \eqref{eq:intro-min-1Segal-square}
  encoding the first Segal condition
  is typically not sent to a Cartesian square by $2$-Segal objects.
  This is explained by \autoref{thm:intro-main-even}:
  while the square \eqref{eq:intro-min-1Segal-square}
  is biCartesian in $\Delta$,
  it is no longer a pushout square in $\Lambda$.
\item
  The $2$-dimensional cyclic polytope $\cycPol{2}{4}$ is a square.
  It has the two triangulations
  (see \autoref{fig:cyc-2-4})
  whose corresponding Segal condition expresses that the two squares
  \begin{equation}
    \label{eq:7}
    \cdsquareNA{13}{123}{013}{0123}
    \intxt{and}
    \cdsquareNA{02}{012}{023}{0123}
  \end{equation}
  in $\Delta$ are sent to a limit diagram.
  Both of the squares \eqref{eq:7} are biCartesian in $\Lambda$.
\item
  The squares
  \begin{equation}
    \label{eq:8}
    \cdsquare
    {11'}{011'}{1}{01}
    {\cofacemap{0}}
    {\codegmap{0}}
    {\codegmap{1}}
    {\cofacemap{0}}
    \intxt{and}
    \cdsquare
    {0'0}{0'01}{0}{01}
    {\cofacemap{2}}
    {\codegmap{0}}
    {\codegmap{0}}
    {\cofacemap{1}}
  \end{equation}
  are biCartesian both in $\Delta$ and in $\Lambda$.
  Hence they need to be sent to pullback squares
  by every Segal object
  (by \autoref{thm:intro-main-odd})
  and by every $2$-Segal object
  (by \autoref{thm:intro-main-even}).
  While the first of these facts is easy, the second is nontrivial;
  it is precisely the statement that $2$-Segal spaces are automatically unital,
  which was discovered only very recently by
  \fivenames{Feller}{Garner}{Kock}{Underhill-Proulx}{Weber}~\cite{FGKUPW2019}.
\end{itemize}

We would like to point out the following
corollary of \autoref{thm:intro-main-even},
which cements the importance of higher cyclic \evenSegal{\kk} objects
and might help explain the particular usefulness of cyclic $2$-Segal objects.

\begin{Corollary}[\autoref{cor:cyclicmain}]
  Let $\localC$ be an $\infty$-category with finite limits.
  The cyclic \evenSegal{\kk} objects in $\localC$
  are precisely the weakly $\kk$-excisive functors $\Lambda^\op\to\localC$.
\end{Corollary}

Finally, we remark that our main theorem implies a nontrivial bound
(\autoref{prop:triviality-bounds})
on how many values of a higher Segal object can be trivial
without the whole object collapsing.
Whether this bound is sharp is still unknown (at least to \theauthor{})
and remains to be investigated in future research.

\subsection{Methods and structure of the paper}

The main conceptual framework which informs our approach
is a version for the simplex category
of the Goodwillie--Weiss manifold calculus.
In \autoref{sec:manifold-calc-Delta} we explain
a system of heuristic analogies between manifold calculus
(in its version described by \nameBdBW{}~\cite{BoavidaWeiss2013})
and a \roughly{manifold calculus} on $\Delta$.
While the mathematics in the rest of the paper stands on its own,
it is \theauthor{}'s opinion that
these informal analogies to manifold calculus
can be very helpful when digesting the definitions and building intuition.
Interestingly, they also explain
how one might have guessed the definition of higher Segal spaces
without knowing about cyclic polytopes.
One practical upshot of the analogy to manifold calculus
is that it inspires the definition of \buzzword{polynomial} simplicial objects,
a notion which is implied by higher weak excision
(while being, \latin{a priori}, weaker)
and which can be compared more easily to the higher Segal conditions.

In \autoref{sec:preliminaries} we recall
basic definitions and facts about
the categories $\Delta$ and $\Lambda$,
(co)Cartesian and strongly (co)Cartesian cubes,
as well as general notions of (weak) excision and descent.
In \autoref{sec:biCart-in-Delta-Lambda},
we explicitly classify strongly Cartesian
and biCartesian cubes in $\Delta$ and in $\Lambda$.
In \autoref{sec:prec-inters-cubes} we explain a
descent theory on $\Delta$ and
and study polynomial simplicial objects in this framework\footnote{
  This framework has already proven its worth
  in the classification of higher Segal objects with values
  in \emph{stable} \inftycats{} by
  \TDyckerhoff{}, \GJasso{} and \theauthor{}~\cite{DJW2018}.
}.
In \autoref{sec:excisive-objects} we show that polynomial simplicial objects
agree with weakly excisive ones;
our key arguments here are a version of the ones in~\cite{FGKUPW2019}
repackaged in a way which directly generalizes to arbitrary dimensions.
The main theorem (\autoref{thm:main})---%
comparing higher Segal conditions with weak excision---%
is proved in the last section (\autoref{sec:higher-Segal})
by considering a series of descent conditions which interpolate
between the higher Segal conditions and the conditions of being polynomial.

\subsection{Acknowledgments}
This work was done during my Ph.D.~studies at the Hausdorff Center for
Mathematics (HCM);
I am very grateful to my supervisors Catharina Stroppel and Tobias Dyckerhoff
for supporting and encouraging me during this time.
This research was supported with a Hausdorff
scholarship by the Bonn International Graduate School for Mathematics (BIGS)
and funded by the Deutsche Forschungsgemeinschaft
(DFG, German Research Foundation)
under Germany's Excellence Strategy - GZ 2047/1, Projekt-ID 390685813.

\subsection{($\infty$-)categorical conventions}
Throughout this article,
we will freely use the theory of $\infty$-categories
as developed in\cite{Lurie2009};
most relevant will be the theory of limits and Kan extensions developed
in chapter 4.
We silently identify each ordinary category with its nerve
so that each ordinary category is in particular an $\infty$-category.
Given two ($\infty$-)categories $\localC$ and $\localD$,
we write $\Fun(\localD,\localC)$ for the ($\infty$-)category of functors
between them;
for instance, $\Fun(\Dop,\localC)$ denotes the $\infty$-category of
simplicial objects in the $\infty$-category $\localC$.
When we talk about commutative diagrams in an $\infty$-category
we will usually only depict objects and arrows
while leaving the higher coherence data implicit.
All limits and colimits are always meant to be taken
in the homotopy-coherent (\ie $\infty$-) sense.

\section{A \roughly{manifold calculus} for the simplex category}
\label{sec:manifold-calc-Delta}

\newcommand{\localMan}{\cat{Man}}
\newcommand{\covertopology}[1]{\calJ_{#1}}
\newcommand{\coveragegood}[1]{\calJ_{#1}^\circ}
\newcommand{\coveragepoly}[1]{\calJ_{#1}^{\mathrm{h}}}

A contravariant functor $\localX$ defined
on the topological (\ie $\infty$-) category $\localMan$
of smooth $\dd$-manifolds and smooth embeddings is usually called
\buzzword{polynomial of degree $\leq 1$}
if its value on a manifold $M$ can be computed
by cutting $M$ up into smaller open pieces,
evaluating $\localX$ piece by piece and then reassembling the values.
More precisely, for each pair of
disjoint closed subsets subsets $A_0,A_1\subset M$,
 one requires the canonical map
\begin{equation*}
  \localX(M)\lra
  \localX(M\setminus A_0)
  \fiberproduct{\localX(M\setminus A_0\cup A_1)}
  \localX(M\setminus A_1)
\end{equation*}
to be an equivalence.

\nameBdBW{}~\cite{BoavidaWeiss2013}
show that polynomial functors of degree $\leq 1$
are precisely the (homotopy) sheaves on $\localMan$
for the Grothendieck topology $\covertopology{1}$ of open covers.
More generally,
they consider a hierarchy $\covertopology{\kk}$
of Grothendieck topologies on $\localMan$
(with $\kk\geq 1$),
where $\covertopology\kk$ consists of those
open covers (called \buzzword{\kk-covers})
which have the property that every set of $\kk$ (or fewer) points is
contained in some open set of the cover.
The manifold calculus of \nameBdBW{}
is concerned with the systematic study of sheaves on
$(\localMan,\covertopology\kk)$.
They introduce the following classes of open
covers:
\begin{enumerate}
\item the class $\coveragepoly\kk$ consists of
  open covers of the form
  \begin{equation}
    \label{eq:manifoldpolycover}
    \setP{M\setminus A_i\hra M}{i=0,\dots,\kk}
  \end{equation}
  for pairwise disjoint closed subsets $A_0,\dots,A_\kk\subset M$ of $M$.
\item the class $\coveragegood\kk$ consists of \buzzword{good} $\kk$-covers,
  \ie $\kk$-covers with the property
  that every finite intersection of open sets
  is diffeomorphic to a disjoint union of at most $\kk$ balls.
\end{enumerate}
While the classes $\coveragepoly\kk$ and $\coveragegood\kk$
are not Grothendieck topologies anymore,
they are so called \buzzword{coverages},
hence they admit a well-behaved theory of descent and sheaves.
Sheaves for the coverage $\coveragepoly\kk$ are called
\buzzword{polynomial functors of degree $\leq \kk$}.
One of the main results of \nameBdBW{} in this context
is the following theorem:

\begin{Thm}\cite[Theorem 5.2 and Theorem 7.2]{BoavidaWeiss2013}
  \label{thm:BWpolygood}
  The coverages $\covertopology\kk$, $\coveragepoly\kk$ and
  $\coveragegood\kk$ define the same class of sheaves on $\localMan$.
\end{Thm}

We shall now describe a similar theory for simplicial objects,
\ie presheaves on the simplex category $\Delta$
(see \autoref{sec:simplex-category} for the notation).
It turns out that the following list of analogies is useful;
we put terms coming from the language of manifold in quotes
to emphasize that they should be thought of heuristically:
\begin{itemize}
\item
  We think of the object $\n=\set{0,\dots,n}\in\Delta$ as a
  \roughly{manifold} with \roughly{points} given by pairs $(x-1,x)$
  with $x=1,\dots,n$.
\item
  An \roughly{open subset} of $\n$ is simply an ordinary subset
  $U\subseteq\set{0,\dots,n}$;
  it contains the \roughly{points} $(x-1,x)$ such that $\set{x-1,x}\subseteq U$.
\item
  We say that two \roughly{open subsets} $U,U'$ of the \roughly{manifold} $\n$
  are \roughly{disjoint} if they are disjoint as subsets of $\n$;
  note that this is a stronger condition than requiring $U$ and $U'$
  to share no \roughly{point}.
\item
  A \roughly{closed set} $A$ of $\n$ is an ordinary subset of $A\subseteq\n$;
  it contains all the points not contained in its complement
  $\n\setminus A\subseteq \n$ (viewed as an \roughly{open set});
  explicitly, $A$ contains all \roughly{points} $(x-1,x)$
  with $x\in A$ or $x-1\in A$.
\item
  We say that two \roughly{closed sets} $A,A'\subseteq\n$
  are \roughly{disjoint} if they share no \roughly{point};
  note that this is stronger than being disjoint as subsets of $\n$.
\item
  Each \roughly{point} $p=(x-1,x)$ has
  a unique minimal \roughly{open neighborhood}
  given by the subset $U^p=\set{x-1,x}\subseteq\n$,
  which we think of as a very small \roughly{open ball}
  around the \roughly{point} $p$.
\end{itemize}

Armed with this intuition, we can define analogs
of the coverings $\coveragepoly\kk$ and $\coveragegood\kk$
in the simplex category:

\begin{enumerate}
\item
  \label{item:7}
  The open covers~\eqref{eq:manifoldpolycover} can be translated to $\Delta$
  by putting everything in quotation marks:
  For every collection $A_0,\dots,A_\kk$
  of \roughly{nonempty and pairwise disjoint closed subsets}
  of the \roughly{manifold} $\n$,
  we can define the \roughly{open cover}
  \begin{equation}
    \label{eq:1}
    \setP{\n\setminus A_i\hra\n}{i=0,\dots,\kk}
  \end{equation}
  of $\n$.
  See also \autoref{sec:polyn-simpl-objects}.
\item
\label{item:14}
  Heuristically\footnote{
    See Proposition 2.10 in~\cite{BoavidaWeiss2013}
    for an actual proof.
  },
  one way to produce good $\kk$-covers of a manifold $M$ is as
  follows: Fix a Riemannian metric on $M$ and, for every tuple
  $p=(p_1,\dots,p_\kk)$ of $\kk$ points in $M$, choose very small (with respect
  to the geodesic distance between the points $p_i$) balls
  $U^p_i\ni p_i$.
  Then the collection
  $\setP{\bigdisjunion_{i=1}^k U^p_i}{p\in M^\kk}$
  is a $\kk$-good cover of $M$.

  In our analogy, every \roughly{point} $p$
  of a \roughly{manifold} $\n\in\Delta$
  has a canonical/minimal \roughly{open ball} $U^p$ surrounding it.
  Hence each $\n\in\Delta$ has a
  canonical \roughly{good $\kk$-cover} containing all those
  \roughly{open subsets} of $\n\in\Delta$
  that can be written as union of the form
  \begin{equation*}
    \bigdisjunion_{i=1}^kU^{p_i},
  \end{equation*}
  where $p_1,\dots,p_k$ are \roughly{points} of the \roughly{manifold} $\n$ with
  \roughly{pairwise disjoint neighborhoods} $U^{p_i}$.
  See also \autoref{sec:higher-segal-covers}.
\end{enumerate}

Inspired by the analogy,
we call a functor $\Dop\to\localC$ \introduce{polynomial of degree $\leq\kk$}
if it is a sheaf for the \roughly{open covers} of type~\ref{item:7}
(see \autoref{def:polynomial-f}).

The following easy observation was \theauthor{}'s original motivation
for this line of inquiry because it shows
on one hand that the canonical \roughly{good $\kk$-covers}
are a meaningful concept
and on the other hand
that a \roughly{manifold calculus} of $\Delta$
can be a powerful organizational principle for higher Segal spaces.

\begin{Obs}
  Sheaves on $\Delta$ with respect to the canonical
  \roughly{good $\kk$-covers} of \ref{item:14}
  are precisely the
  \loweroddSegal{\kk} spaces of \twonames{Dyckerhoff}{Kapranov}.
\end{Obs}

The notion of polynomial simplicial objects might be a bit unsatisfying
because its very definition relies on an informal analogy to manifold calculus;
without this analogy, the \roughly{open covers}~\eqref{eq:manifoldpolycover}
might seem a bit mysterious and devoid of intrinsic meaning.
We will clarify this issue by showing that a functor $\Dop\to\localC$
is polynomial of degree $\leq\kk$
if and only if it is weakly $\kk$-excisive
(see \autoref{thm:unitality-of-cubes}).
In this light, our main result (\autoref{thm:main})
relating \loweroddSegal{\kk} objects with weakly $\kk$-excisive functors
should be seen as a discrete analog of \autoref{thm:BWpolygood}
of \twonames{Boavida de Brito}{Weiss}.\\

We will not spell out the whole story for \evenSegal{\kk} objects
since it is very similar.
Let us just say that one should now consider a \roughly{manifold calculus}
not on the simplex category $\Delta$ but on \ConnesCC{} $\Lambda$,
where the \roughly{manifold} $\n=\set{0,\dots,n}$
now has an additional \roughly{point} given by $(n,0)$.

\section{Preliminaries}
\label{sec:preliminaries}

\subsection{The simplex category}
\label{sec:simplex-category}
The \introduce{augmented simplex category} $\Deltaug$ is the category
of finite linearly ordered sets and order preserving
(\ie weakly monotone)
maps between them.
The \introduce{simplex category} $\Delta\subset\Deltaug$ is the full subcategory
spanned by the \emph{nonempty} finite linearly ordered sets.
Every object in $\Delta$ is isomorphic, by a unique isomorphism,
to a standard ordinal of the form
$\n\coloneqq\set{0<1<\dots< n}$ for some $n\in\BN$;
when convenient can we therefore identify $\Delta$
with its skeleton spanned by $\setP{\n}{n\in\BN}$.

\begin{Def}
  A \introduce{simplicial object} in an ($\infty$-)category $\localC$
  is a functor $\Dop\to\localC$.
\end{Def}

The augmented simplex category has a monoidal structure
\begin{equation}
  \label{eq:joinonDp}
  \join\colon\Deltaug\times\Deltaug\lra\Deltaug,
\end{equation}
given by left-to-right concatenation or \introduce{join}
of linearly ordered sets.
Explicitly we have
\begin{equation*}
  \set{a_0<\dots<a_n}\join \set{b_0<\dots<b_m}
  \coloneqq\set{a_0<\dots<a_n<b_0<\dots<b_m};
\end{equation*}
the monoidal unit for $\join$ is the empty set $\emptyset\in\Deltaug$.
We use the convention
$\numD{-1}\coloneqq\emptyset\in\Deltaug$ and
$\numDminus{n}{m}\coloneqq\set{m+1<\dots< n}$
for all $-1\leq m\leq n$
so that we always have
$\n=\m\join\numDminus nm$.
Given a simplicial object $\xdc$,
the \introduce{left path object $\lpath\localX$}
and the
the \introduce{right path object $\rpath\localX$}
are defined as the compositions
\begin{equation*}
  \lpath\localX\colon
  \Dop\xra{\0\join\blank}
  \Dop\xra{\localX}\localC
  \intxt{and}
  \rpath\localX\colon
  \Dop\xra{\blank\join\0}
  \Dop\xra{\localX}\localC,
\end{equation*}
respectively.

A morphism $f\colon\m\to\n$ in $\Delta$ is called
\introduce{left active} or \introduce{right active},
if it preserves the minimal element (\ie $f(0)=0$)
or the maximal element (\ie $f(m)=n$), respectively;
call $f$ \introduce{active} if it is both left and right active.
Denote by $\Dlact$, $\Dract$ and $\Dact\coloneqq\Dlact\cap\Dract$
the wide subcategories of $\Delta$ containing
the left active, right active and active morphisms, respectively.
Call a morphism $f\colon \m\to\n$
\introduce{left strict} or \introduce{right strict}
if $\preimage f0=\set0$ or $\preimage fn= \set m$, respectively.
For each $n\in\BN$,
we denote by $\canactivemap n\colon \numD 1\to\n$ the unique active map;
explicitly given as $\canactivemap n(0)=0$ and $\canactivemap{n}(1)=n$.

\subsection{The cyclic category}
\newcommand{\excycset}{N}
\newcommand{\excycendo}{T}
\newcommand{\exlinsubset}{L}
\newcommand{\exunderlinsubset}{L_0}
\newcommand{\linpreimage}[2]{\lininvimage{#1}\set{#2}}
\newcommand{\lininvimage}[1]{{#1}^\star}
\newcommand{\cff}f%

  A \introduce{finite cyclic set} is a pair $(\excycset,\excycendo)$
  consisting of a finite set $\excycset$
  together with an endomorphism
  $\excycendo\colon\excycset\to\excycset$ which is
  transitive, \ie for each $x,y\in\excycset$ there is some $i\in\BN$ such that
  $\excycendo^ix=y$.
  A \introduce{linearly ordered subset}
  $\exlinsubset=(\exunderlinsubset,\prec)$ of
  $(\excycset,\excycendo)$ is a subset $\exunderlinsubset$ of $\excycset$
  (called the \introduce{underlying set} of $\exlinsubset$)
  equipped with a linear order $\prec$ such that
  the restriction of $\excycendo$ to $\exlinsubset$
  agrees with the successor function induced by $\prec$.
  A morphism
  $(\cff,\lininvimage{\cff})\colon
  (\excycset',\excycendo')\lra(\excycset',\excycendo')$
  of finite cyclic sets consists of
  \begin{itemize}
  \item
    a map of sets $\excycset'\to\excycset$ which we also denote by $\cff$
    and
  \item
    an assignment $\lininvimage\cff$, which
    for each linearly ordered subset $\exlinsubset\subset \excycset$ produces
    a linearly ordered subset $\lininvimage\cff \exlinsubset\subset \excycset'$
    with underlying set $\invimage\cff\exlinsubset$ such that
    $\lininvimage\cff{\exlinsubset}
    =\lininvimage\cff{\exlinsubset'}\join \lininvimage\cff{\exlinsubset''}$
    whenever the linerly ordered subset $\exlinsubset\subset\excycset$ is
    decomposed as $\exlinsubset=\exlinsubset'\join\exlinsubset''$.
  \end{itemize}
  Composition of morphisms
  $\excycset''
  \xra{(\cff',\lininvimage{\cff'})}\excycset'
  \xra{(\cff,\lininvimage{\cff})}\excycset$
  between finite cyclic set is given
  by the usual composition of underlying set maps and
  $\lininvimage{(\cff\circ\cff')}=\lininvimage{\cff'}\circ\lininvimage\cff$.

\begin{Def}[\cite{Connes1983}]
  \introduce{\ConnesCC{}}
  $\Lambda$
  consists of nonempty finite cyclic sets and morphisms between them.
  A \introduce{cyclic object} in some ($\infty$-)category $\localC$
  is a functor $\localX\colon\Lambda^\op\to\localC$.
\end{Def}

For each $n\in \BN$, we have the standard finite cyclic set
\begin{equation*}
  \numL n\coloneqq\left(\cygrp{(n+1)},{+1}\right).
\end{equation*}
It is easy to see that every nonempty finite cyclic set is (non-canonically)
isomorphic to exactly one such standard cyclic set.
Motivated by this, we use the notation ${+m}\coloneqq \excycendo^m$ and
${-1}\coloneqq\excycendo^{-m}$ for arbitrary finite cyclic sets
$(\excycset,\excycendo)$ and omit $\excycendo$ from the notation entirely.

For every finite cyclic set $(\excycset, +1)$, the automorphism group
$\Aut_\Lambda(\excycset)$ is cyclic of order $\card{\excycset}$
and is generated by the structure morphism $+1\colon\excycset \to\excycset$
where $\lininvimage{(+1)}\coloneqq {-1}$ is given by
\begin{equation*}
  \excycset\supset\exlinsubset\lmapsto\exlinsubset -1\coloneqq
  \setP{x-1}{x\in \exlinsubset}\subset\excycset.
\end{equation*}

Specifying a morphism $\cff\colon\excycset\to\numL0$ amounts to the choice of
what we call a \introduce{linear order on the cyclic set} $\excycset$,
namely a linearly ordered subset $\linpreimage\cff 0\subset\excycset$
with underlying set $\preimage\cff 0=\excycset$.
A commutative triangle
\begin{equation*}
  \begin{tikzcd}
    \excycset'\ar[rd,"\cff'"']\ar[rr]&&\excycset\ar[dl,"\cff"]\\
    &\numL 0
  \end{tikzcd}
\end{equation*}
corresponds precisely to an order preserving map
$\linpreimage{\cff'} 0\to \linpreimage\cff 0$.
We conclude that the assignment
$\cff\mapsto \linpreimage{\cff}0$ describes a functor
\begin{equation*}
  \overcat\Lambda {\numL 0}\xra{\cong}\Delta,
\end{equation*}
which is easily seen to be an isomorphism of categories
between $\Delta$ and the slice of $\Lambda$ over $\numL 0$.
\newcommand{\structuremapnumL}[1]{\numD{#1}}
Under this identification,
the object $\n\in\Delta$ corresponds to $\numL n\in\Lambda$
which is equipped with the structure map
$\structuremapnumL n\colon\numL n\to \numL 0$
induced by the standard linear order $0<1<\dots< n$ on $\cygrp{(n+1)}$.

\newcommand{\cyclicrotation}[2]{{#2}^{#1}}
\newcommand{\Dcyclicrotation}[2]{\structuremapnumL{#2}^{#1}}

Composition in $\Lambda$ induces a free and transitive right group action
\begin{align*}
  \label{eq:autLambdaaction}
  \Lambda(\excycset,\numL 0)\times \Aut_\Lambda(\numL n)
  &\lra\Lambda(\excycset,\numL 0);\\
  (\cff,{+m})&\lmapsto \cyclicrotation{+m}{\cff}
\end{align*}
which corresponds to cyclic rotation of linear orders:
if $\structuremapnumL n\colon\numL n\to\numL0$ is the structure map
corresponding to the standard order $<$ on $\n$,
then $\Dcyclicrotation{+m}n$ corresponds to the linear order $\prec$ on the set
$\set{0,1\dots, n}$ given by
${n-m+1}\prec \dots \prec n\prec 0\prec \dots\prec {n-m}$.

\subsection{Cartesian and coCartesian cubes}
\label{sec:cart-cocart-cubes}

Fix a finite set $S$ and denote by $\powerset S$ the powerset of $S$, partially
ordered by inclusion.

\begin{Def} An \introduce{$S$-cube} in some ($\infty$-)category $\localC$
  is a functor $\excube\colon\powerset S\to \localC$.
\end{Def}

\begin{Rem}
  Since the poset $\powerset S$ is canonically isomorphic to its opposite
  (via the assignment $S\supseteq T\mapsto {S\setminus T}$),
  we will often write cubes in an ($\infty$-)category $\localD$
  as functors $\powersetop S\to \localD$.
  This is convenient when studying contravariant functors
  $\localX\colon\localD^\op\to\localC$,
  where we can then say that the cube
  $\powersetop S\to\localD$ in $\localD$ is sent by $\localX$ to the composite
  $\powerset S\to \localD^\op\xra{\localX}\localC$;
  the main example in this paper is of course the case where $\localD=\Delta$
  and $\xdc$ is a simplicial object in $\localC$.
\end{Rem}

\newcommand{\cubeindirection}[2]{#1^{#2}}
\newcommand{\cubeindirectionS}[2]{#1^{#2\notin}}
\newcommand{\cubeindirectionT}[2]{#1^{#2\in}}

Let $s\in S$ and write $S'\coloneqq S\setminus {\set s}$.
For every element $s\in S$ we have an isomorphism of posets
\begin{equation}
  \label{eq:sSisocube}
  {\simplex 1}\times\powerset{S'}\xra{\cong} \powerset{S}
\end{equation}
given by $(0,T)\mapsto T$ and $(1,T)\mapsto T\disjunion{\set s}$.
For every $\infty$-category $\localC$ we get an induced equivalence
\begin{equation}
  \label{eq:unravelcube}%
  \Fun(\powerset S,\localC)\xra{\simeq}%
  \Fun(\simplex 1,\Fun(\powerset{S'},\localC))%
\end{equation}
of $\infty$-categories,
which we denote by $\excube\mapsto \cubeindirection{\excube}{s}$.
We say that a cube $\excube$ is the pasting in $s$-direction of two cubes
$\excube'$ and $\excube'$ if we have an identification
$\cubeindirection{\excube}{s}=
\cubeindirection{\excube'}{s}\circ
\cubeindirection{\excube''}{s}$.

Denote by $\powersetne S\coloneqq\powerset S\setminus{\set\emptyset}$
the poset of \emph{nonempty} subsets of $S$.

\begin{Def}
  \label{def:coCartcubes}
  An $S$-cube $\excube\colon\powerset S\to \localC$ is called
  \begin{itemize}
  \item
    \introduce{Cartesian} if it is a limit diagram in $\localC$;
    \ie if $\excube$ is the right Kan extension
    of its restriction to $\powersetne S$.
  \item
    \introduce{coCartesian} if it is a colimit diagram in $\localC$;
    \ie if $\excube$ is the left Kan extension
    of its restriction to $\powerset S\setminus{\set S}$.
    \qedhere
  \end{itemize}
\end{Def}

\begin{Def}
  \label{def:strongly-biCart}
  An $S$-cube $\excube\colon\powersetop S\to\localD$ is called
  \introduce{strongly Cartesian} or \introduce{strongly coCartesian}
  if,
  for each $T\subset S$ and $s,s'\in S\setminus T$ with $s\neq s'$,
  the $2$-dimensional face
  \begin{equation}
    \label{eq:faceforcoCart}
    \begin{tikzcd}
      T\ar[r]\ar[d]&T\disjunion\set s\ar[d]\\
      T\disjunion\set{s'}\ar[r]&T\disjunion\set{s,s'}
    \end{tikzcd}
  \end{equation}
  is sent by $\excube$ to
  a pullback square or a pushout square in $\localD$,
  respectively.
  A cube is called \introduce{strongly biCartesian} if it is both
  strongly Cartesian and strongly coCartesian.
\end{Def}

\begin{Rem}
  Denote by $\powersetleop S1$ and by $\powersetgeop S{\card S-1}$
  the subposet of $\powersetop S$ spanned by the
  subsets $T\subset S$ of cardinality
  $\card T\leq 1$ and $\card T\geq \card S-1$, respectively.
  It is easy to see that a cube $\excube\colon \powersetop S\to\localD$
  is strongly Cartesian if and only if it is the right Kan extension of its
  restriction to $\powersetleop S1$;
  it is strongly coCartesian if and only if it is the left Kan extension of its
  restriction to $\powersetgeop S{\card S-1}$.
\end{Rem}

\begin{Rem}
  If $\card S\geq 2$, then every strongly (co)Cartesian cube is also
  (co)Cartesian; thus justifying the terminology. Beware however, that for
  $\card S=1$ an $S$-cube is just an arrow;
  it is always strongly biCartesian and
  is (co)Cartesian if and only if it is an equivalence.
\end{Rem}

\begin{Lem}
  \label{lem:cartesiancubesasedges}
  Let $\localC$ be an $\infty$-category. Let $s\in S$ and put $S'\coloneqq
  S\setminus \set s$. The restriction functor
  \begin{equation}
    \label{eq:respuncturedcube}
    p\colon\Fun(\powerset {S'}, \localC)\lra \Fun(\powersetne {S'},\localC)
  \end{equation}
  is a coCartesian fibration which is Cartesian if and only if $\localC$ admits
  limits of shape $\powersetne S$.
  An $S$-cube $\excube\colon\powerset S\to\localC$ is Cartesian if and only if
  the corresponding edge
  $\cubeindirection\excube s\colon\simplex1\to\Fun(\powerset{S'},\localC)$
  is $p$-Cartesian.
\end{Lem}

\begin{Prf}
  \autoref{lem:cartesiancubesasedges} is the higher dimensional analog of
  Lemma~6.1.1.1 in~\cite{Lurie2009};
  the proof is essentially the same.
\end{Prf}

We say that an  $S$-cube $\excube$ is
\introduce{degenerate in direction $s\in S$}
if the corresponding natural transformation $\cubeindirection\excube s$
of $S\setminus {\set s}$-cubes is an equivalence.
It follows directly from~\autoref{lem:cartesiancubesasedges} that
\introduce{degenerate cubes}---cubes
that are degenerate in at least one direction---are
automatically Cartesian and coCartesian.

The following lemma is a standard argument which is useful to compare
Cartesian cubes of different dimensions.

\begin{Lem}
  \label{lem:changeofclaw}
  Let $\excube\colon\powerset S\to\localC$ be an $S$-cube in
  an $\infty$-category $\localC$ with finite limits.
  Fix $s\in S$ and write $S'\coloneqq {S\setminus{\set s}}$.
  Assume that the $S'$-cube
  $\cubeindirection\excube s(1)\colon T\mapsto \excube(T\disjunion {\set s})$
  is Cartesian.
  Then the canonical map
  \begin{equation}
    \label{eq:canmapcuberestriction}
    \lim \restr\excube{\powersetne S}
    \lra
    \lim\restr{\excube}{\powersetne{S'}}
  \end{equation}
  is an equivalence.
  In particular, the original $S$-cube $\excube$ is Cartesian
  if and only if the restricted $S'$-cube
  $\restr\excube{\powerset{S'}}
  =\cubeindirection\excube s(0)
  \colon T\mapsto \excube(T)$
  is Cartesian.
\end{Lem}

\begin{Prf}
  Consider the following commutative diagram in $\localC$
  \begin{equation}
    \label{eq:diainproofnewdirection}
    \begin{tikzcd}
      Q(\emptyset)\ar[r]\ar[rd]&
      \lim \restr{\excube}{\powersetne S}\ar[r]\ar[d]\isCartesian&
      \excube(\set s)\ar[d,"{\simeq}"]\\
      &\lim \restr{\excube'}{\powersetne{S'}}\ar[r]&
      \lim \restr{\cubeindirection\excube s(1)}{\powersetne {S'}}
    \end{tikzcd}
  \end{equation}
  which is induced by the universal properties of the various limits.
  By a standard decomposition argument for limits%
  , the rightmost square in the diagram~\eqref{eq:diainproofnewdirection}
  is Cartesian;
  moreover, the rightmost vertical map is an equivalence by assumption.
  It follows that the left vertical map is also an equivalence;
  the result follows.
\end{Prf}

\subsection{\Cech{} cubes, descent and weak excision}
\label{sec:weak-excision}

\newcommand{\exNinD}{N}
Let $\localD$ be an \inftycat{}.

\begin{Def}
  Let $S$ be a finite set.
  An \introduce{$S$-pronged claw} (or just \introduce{$S$-claw}, for short)
  $\exclaw$ on an object $\exNinD$ in $\localD$
  is an $S$-indexed tuple
  $\exclaw=\tupleP{\ff_s\colon\II_s\to\exNinD}{s\in S}$
  of maps $\ff_s$ in $\localD$ with common codomain $\exNinD\in\localD$
  or, equivalently,
  a diagram $\exclaw\colon\powersetleop S1\to \localD$
  with $\exclaw(\emptyset)=\exNinD$.
\end{Def}

Given an $S$-claw
$\exclaw=\tupleP{\ff_s\colon\II_s\to \exNinD}{s\in S}$ on $\exNinD\in\localD$,
we write $\exclaw\clawon \exNinD$
to make the codomain $\exNinD$ explicit in the notation
while keeping the $\ff_s$, the $\II_s$ and sometimes even the $S$ anonymous.
In a similar spirit we will use the symbol
$\ff\in\exclaw$ to mean $\ff_s$ for some $s$.
With this convention $\ff_s$ and $\ff_{s'}$ should be considered distinct if
$s\neq s'$, even if they are the same map in $\localD$.
Each subset $T\subset S$ induces a restricted $T$-claw of $\exclaw$ given by
$\restr{\exclaw}{T}\coloneqq\tupleP{\ff_t}{t\in T}\clawon\exNinD$.

\begin{Def}
  An $S$-claw $\exclaw\clawon\exNinD$ in $\localD$ is called
  a \introduce{candidate $S$-covering}
  if it can be extended to a strongly Cartesian $S$-cube
  $\Cechcube\exclaw\colon\powersetop{S}\to\localD$.
  In this case we call $\Cechcube\exclaw$
  the \introduce{\Cech{} cube} associated to $\exclaw$.
\end{Def}

If it exists, the \Cech{} cube $\Cechcube\exclaw$ is given by the formula
\begin{equation}
  \label{eq:formulaCech}
  S \supseteq T\lmapsto\lim{\restr{\exclaw}{T}}.
\end{equation}
We shall sometimes think of the prongs $\ff_s\colon\II_s\to\exNinD$
as generalized subobjects of $\exNinD$;
the values \eqref{eq:formulaCech} of the \Cech{} cube should then be thought of
as generalized intersections.
In this spirit it is sometimes convenient to use the notation
$\bigcap_{t\in T}\ff_t\coloneqq \Cechcube\exclaw(T)=\lim{\restr{\exclaw}{T}}$
and denote, for instance, the \Cech{} square of two maps
$\ff\colon\II\to\n$ and $\ff'\colon\II'\to\n$ as follows:
\begin{equation}
  \label{eq:abusiveintersection}
  \begin{tikzcd}
    \II\cap\II'\rar{\ff\cap\II'}\ar[d, "\II\cap\ff'"']
    \ar[dr,"\ff\cap\ff'" description]
    & \II'\dar{\ff'}\\
    \II\rar{\ff}&\n
  \end{tikzcd}
\end{equation}

\begin{Def}
  Let $\exclaw$ be a candidate covering in $\localD$.
  A functor $\localX\colon\localD\to\localC$ is said to
  \introduce{satisfy descent with respect $\exclaw$}
  if it sends the \Cech{} cube $\Cechcube\exclaw$
  to a Cartesian cube in $\localC$;
  in this case we also say that $\exclaw$ is \introduce{$\localX$-local}.
\end{Def}

\newcommand{\excoverage}{\tau}

Following \nameBdBW{} we say
that a \introduce{coverage} $\excoverage$ on $\localD$ is a
collection of candidate coverings.
If $\exclaw\clawon\exNinD$ is an element of $\excoverage$
then we say that $\exclaw$ is a \introduce{$\excoverage$-covering};
if the coverage $\excoverage$ is implicit from the context
then we say that $\exclaw$ is a \introduce{covering} of $\exNinD$.

\begin{Def}
  A $\localC$-valued \introduce{sheaf}
  for the coverage $\excoverage$ is a functor
  $\localX\colon\localD^\op\to \localC$
  which satisfies descent with respect to all $\excoverage$-coverings.
\end{Def}

\newcommand\coverageall[1]{\excoverage_{#1}}

\begin{Rem}
  For each $\kk\geq 0$,
  there is a canonical coverage $\coverageall \kk$ on $\localD$
  which consists of \emph{all} candidate $\numD\kk$-coverings.
  A presheaf $\localD^\op\to\localC$ is
  a sheaf for this coverage $\coverageall \kk$
  if and only if
  it is an $\kk$-excisive (covariant) functor
  in the sense of Goodwillie~\cite{Goodwillie91},
  \ie if it sends strongly coCartesian $\numD\kk$-cubes in $\localD^\op$
  to Cartesian cubes in $\localC$.
\end{Rem}

We say that an $S$-claw is
\introduce{strongly biCartesian} if it is a candidate covering and
if its \Cech{} cube is strongly coCartesian
(hence strongly biCartesian).

\begin{Def}
  A functor $\localD^\op\to\localC$ is called \introduce{weakly $S$-excisive}
  if it is a sheaf for the coverage of strongly biCartesian $S$-claws,
  \ie if it sends all strongly biCartesian $S$-cubes
  to Cartesian cubes in $\localC$.
\end{Def}

We will also need the following relative notion:

\begin{Def}
  \label{def:relative-exc}
  Let $\localD\to\localD'$ be a limit preserving functor.
  We call a functor $\localX\colon\localD^\op\to\localC$
  \introduce{weakly $S$-$\localD'$-excisive}
  (with the functor $\localD\to\localD'$ left implicit)
  if it is a sheaf with respect to those candidate $S$-coverings
  which become strongly biCartesian in $\localD'$.
\end{Def}

Clearly the property of being weakly $S$-excisive
(both in the relative and in the absolute sense)
only depends on the cardinality of $S$.
For $\kk\in \BN$, we say that $\xdc$ is \introduce{weakly $\kk$-excisive}
if it is weakly $\numD{\kk}$-excisive.
We will stick to $S$-cubes instead of $\numD\kk$-cubes whenever possible,
because the latter might suggest
a dependency on the linear order of the coordinates.

\begin{Rem}
\label{rem:restr-relative-exc}
  In the setting of \autoref{def:relative-exc},
  if every candidate covering in $\localD'$ admits
  a lift to a candidate covering in $\localD$
  then a functor $\localD'^\op\to\localC$
  is weakly $S$-excisive if and only if its restriction to $\localD$
  is weakly $S$-$\localD'$-excisive.
\end{Rem}

\section{Strongly biCartesian cubes in $\Delta$ and $\Lambda$}
\label{sec:biCart-in-Delta-Lambda}

The goal of this section is to classify and explicitly describe the strongly
biCartesian cubes in the simplex category and the cyclic category.

\subsection{Strongly biCartesian cubes in the simplex category}

\begin{Def}
  An $S$-claw $\exclaw=\tupleP{\ff_s}{s\in S}$ on $\n$ in $\Deltaug$
  is called
  \begin{itemize}
  \item%
    \introduce{backwards compatible} if
    for each $i\in\n$ there is at most one $s\in S$ such that the preimage
    $\preimage{\ff_s}{i}$ has more than one element;
  \item
    \introduce{compatible} if it satisfies the following two conditions:
    \begin{enumerate}[label=\axiomlabelstyle{BC}]
    \item
      \label{ax:preimagesforbiCart}
      for each $i\in\n$, there is at most one $s\in S$
      such that the preimage $\preimage{\ff_s}{i}$ is not a singleton;
    \item
      \label{ax:gapsforpushout}
      for each $0<i\leq n$, there is at most one $s\in S$
      such that the subset $\set{i-1,i}\subseteq\n$
      is not contained in the image of $\ff_s$.
      \qedhere
    \end{enumerate}
  \end{itemize}
\end{Def}

\begin{Rem}
  The $S$-claw $\exclaw$ satisfies condition~\ref{ax:preimagesforbiCart}
  if and only if it is backwards compatible and:
  if the preimage $\preimage{\ff_s}i$ is empty for some $i\in \n$ and $s\in S$
  then the preimage $\preimage{\ff_{s'}}i$ is a singleton
  for all $s'\in S\setminus{s}$.
  In the language of \autoref{sec:manifold-calc-Delta},
  condition~\ref{ax:gapsforpushout} says precisely
  that the images of the maps $\ff_s$ are of the form $\n\setminus A_s$,
  where the $\tupleP{A_s}{s\in S}$ are
  \roughly{pairwise disjoint closed subsets} of the \roughly{manifold} $\n$.
\end{Rem}

We call a diagram in $\Deltaug$
\introduce{left active} or \introduce{right active}
if it takes values in the subcategory of $\Delta$
spanned by the left active or right active morphisms, respectively.

\newcommand{\xxX}{*}
\newcommand{\xxO}{\emptyset}
\newcommand{\xx}[1]{#1}

\begin{Rem}
  \label{rem:visualization} 
  It will be useful to visualize $S$-claws
  $\exclaw\clawon\n$ as arrays as in the following example
  (with $n=9$ and $S=\numD 3$):
  \begin{equation}
    \label{eq:visualization-compatible}
    \begin{array}{c|*{10}c}
        & 0 & 1 & 2 & 3 & 4 & 5   & 6  & 7 & 8  & 9\\
      \hline
      0 &\xxX&\xxX&\xxX&\xxX&\xxX&\xx3&\xx2&\xxX&\xxX&\xxX\\
      1 &\xxX&\xxO&\xxX&\xxX&\xxX&\xxX&\xxX&\xxX&\xxX&\xxX\\
      2 &\xxX&\xxX&\xxX&\xxO&\xxO&\xxX&\xxX&\xxX&\xxX&\xxX\\
      3 &\xxX&\xxX&\xxX&\xxX&\xxX&\xxX&\xxX&\xxO&\xx2&\xxO
    \end{array}
  \end{equation}
  There is one row for each prong $\ff_s\colon\II_s\to\n$
  of $\exclaw$ and one column for each $i\in\n$;
  in the column $(s,i)$ we draw:
  \begin{itemize}
  \item
    a star $\xxX$ if the preimage $\preimage{\ff_s}{i}$ is a singleton,
  \item
    the symbol $\xxO$ if the preimage $\preimage{\ff_s}{i}$ is empty or
  \item
    a number $l$ if the preimage $\preimage{\ff_s}{i}$ has $l>1$ many elements.
  \end{itemize}
  A claw is backwards compatible if and only if
  in each column there is at most one entry
  with more than one star.
  It is compatible if and only if
  it satisfies the following two conditions:
  \begin{itemize}
  \item
    in each column there is at most one \roughly{special} entry,
    \ie a cell which is not a star $*$;
  \item
    each pair of two empty cells is either in the same row
    or separated by a column with no empty cells.
  \end{itemize}
  The example \eqref{eq:visualization-compatible}
  depicts a left active compatible claw.
\end{Rem}

\begin{Prop}
  \label{prop:coCartCechcubes}
  Let $\exclaw\clawon\n$ be an $S$-claw in $\Deltaug$.
  \begin{enumerate}[label=(\alph*)]
  \item\label{item:Cechcube}
    The claw $\exclaw$ is a candidate $S$-covering in $\Deltaug$
    if and only if $\exclaw$ is backwards compatible.
    The \Cech{} cube $\Cechcube\exclaw\colon\powersetop S\to\Deltaug$
    is given explicitly by the formula
    \begin{equation}
      \label{eq:explicitpullback}
      \Cechcube\exclaw\colon T\lmapsto
      \bigjoinl_{i\in\n}\prod_{t\in T}\preimage{f_t}{i}.
    \end{equation}
  \item
    The $S$-claw $\exclaw$ isstrongly biCartesian
    (\ie the \Cech{} cube $\Cechcube\exclaw$
    of $\exclaw$ is strongly biCartesian)
    if and only if $\exclaw$ is compatible.
    \qedhere
  \end{enumerate}
\end{Prop}

\begin{Cor}
  \label{cor:biCart=comp}
  A claw in $\Delta$ is strongly biCartesian if and only if it is compatible.
\end{Cor}

\begin{Prf}
  \autoref{cor:biCart=comp} follows directly from \autoref{prop:coCartCechcubes}
  and the easy observation that
  the whole \Cech{} cube of a \emph{compatible} claw
  $\exclaw\clawon\n$ in $\Deltaug$
  lies in $\Delta$ provided that $n\neq -1$.
\end{Prf}

\begin{Expl}
  \label{expl:1-Segal-claw}
   The lower $\numD 1$-claw
  \begin{equation}
    \label{eq:2}
    \begin{array}{*{3}c}
      0 & 1 & 2\\
      \hline
      \xxO&\xxX&\xxX\\
      \xxX&\xxX&\xxO
    \end{array}
  \end{equation}
  is compatible and gives rise to the biCartesian square
  \begin{equation}
    \label{eq:3}
    \cdsquare[bC]{1}{12}{01}{012}{}{}{}{}
  \end{equation}
  in $\Delta$ which encodes the lowest instance of Rezk's Segal conditions.
\end{Expl}

\begin{Prf}[of~\autoref{prop:coCartCechcubes}]
  \begin{enumerate}[label=(\alph*)]
  \item
    \label{item:arisentrs}
    \latin{A priori}, the formula~\eqref{eq:explicitpullback}
    describes a strongly Cartesian extension
    $\Cechcube\exclaw\colon\powersetop S\to \poSet$ of $\exclaw$
    in the category of posets.
    Since the canonical inclusion $\Deltaug\hra\poSet$
    preserves limits, we conclude that $\Cechcube\exclaw$ is
    a strongly Cartesian extension of $\exclaw$ in $\Deltaug$
    if and only if $\Cechcube\exclaw$ takes values
    in linearly ordered posets. This happens if and only if each product
    $\prod_{t\in T}\preimage{\ff_t}{i}$
    has at most one factor which is not empty or a singleton;
    this is precisely the backwards compatibility condition on $\exclaw$.
  \item
    Assume that $\exclaw$ is backwards compatible so that the \Cech{} cube
    $\Cechcube\exclaw\coloneqq\powersetop S\to \Deltaug$ is well defined
    by part~\ref{item:arisentrs}.
    We need to understand when $\Cechcube\exclaw$
    is additionally strongly coCartesian.
    By definition, the cube $\Cechcube\exclaw$ is strongly coCartesian
    if and only it for every subset $T\subset S$ and every pair of
    distinct elements $s,s'\in S\setminus T$, the square
    \begin{equation}
      \label{eq:tocheckstrcocart}
      \begin{tikzcd}[column sep=small]
        \bigjoinl_{i\in\n}%
        \left(\preimage{\ff_s}{i}\times%
          \preimage{\ff_{s'}}{i}\times%
          \prod_{t\in T}\preimage{\ff_t}{i}\right)%
        \ar[r]\ar[d]&%
        \bigjoinl_{i\in\n}%
        \left(\preimage{\ff_{s'}}{i}\times%
          \prod_{t\in T}\preimage{\ff_t}{i}\right)\ar[d]\eqqcolon B'\\%
        B\coloneqq\bigjoinl_{i\in\n}%
        \left(\preimage{\ff_{s}}{i}\times%
          \prod_{t\in T}\preimage{\ff_t}{i}\right)\ar[r]&%
        \bigjoinl_{i\in\n}%
        \left(\prod_{t\in T}\preimage{\ff_t}{i}\right)\eqqcolon N%
      \end{tikzcd}
    \end{equation}
    is a pushout in $\Deltaug$.

    To show \qquote{if} in the claimed equivalence,
    assume that $\exclaw$ is compatible; we will show that then each
    square~\eqref{eq:tocheckstrcocart} is a pushout in $\Deltaug$.
    Condition~\ref{ax:preimagesforbiCart} implies that, for every $i\in\n$,
    if one amongst $\preimage{\ff_s}i$ and $\preimage{\ff_{s'}}i$ is empty
    then the other is a singleton; it follows that the
    square~\eqref{eq:tocheckstrcocart}
    is a pushout on the level of underlying sets.
    It remains to show that a map of sets $\beta\colon N\to M$
    is weakly monotone if it is weakly monotone when
    composed with $B\to N$ and $B'\to N$;
    for this it is sufficient to show that each pair of adjacent elements
    in $N$ is contained in the image of $B\to N$ or in the image of $B'\to N$.
    Let $x<x+1\eqqcolon x'$ be two adjacent elements of $N$
    and denote by $i$ and $i'$ their respective images in $\n$.
    It is enough to show that the subset $\set{i,i'}\subseteq\n$
    is contained in the image of $\ff_s$ or in the image of $\ff_{s'}$.
    If $i=i'$ then this follows from condition~\ref{ax:preimagesforbiCart};
    if $i'=i+1$ then this follows from condition~\ref{ax:gapsforpushout}.
    We may therefore assume $i<i+1\leq i'-1<i'$.
    For each $i<i''<i'$ the product
    $\prod_{t\in T}\preimage{\ff_t}{i''}$
    must be empty by adjacency of $x$ and $x'$.
    Hence there must be $t,t'\in T$ such that
    $\preimage{\ff_t}{i+1}$ and $\preimage{\ff_{t'}}{i'-1}$
    are empty;
    in particular the subsets $\set{i,i+1}$ and $\set{i'-1,i}$ of $\n$ are not
    contained in the image of $\ff_t$ and $\ff_{t'}$, respectively.
    Condition~\ref{ax:gapsforpushout} implies that
    the sets $\set{i,i+1}$, $\set{i'-1,i}$ and, \latin{a fortiori},
    $\set{i,i'}$ are contained in the image of both $\ff_s$ and $\ff_{s'}$.

    To show \qquote{only if}, assume that
    the cube $\Cechcube\exclaw$ is strongly biCartesian. We show that
    conditions \ref{ax:preimagesforbiCart} and \ref{ax:gapsforpushout} hold,
    \ie that $\exclaw$ is compatible.
    \begin{itemize}
    \item[\ref{ax:preimagesforbiCart}]
      Let $i\in\n$ and $s\in S$ be such that $\preimage{\ff_s}{i}$ is empty.
      For each $s'\in S\setminus{\set s}$
      consider the following commutative diagram,
      where the inner solid square is the pushout
      square~\eqref{eq:tocheckstrcocart} (for $T=\emptyset$):
      \begin{equation*}
        \begin{tikzcd}[column sep=small]
          \bigjoinl_{j\in\n}%
          \preimage{\ff_s}{j}\times%
          \preimage{\ff_{s'}}{j}%
          \ar[r]\ar[d]&%
          \bigjoinl_{j\in\n} \preimage{\ff_{s'}}{j}%
          \ar[d]\ar[rdd]\\%
          \bigjoinl_{j\in\n} \preimage{\ff_{s}}{j}%
          \ar[r]\ar[d,equal]&%
          \n\ar[dr,dashed]\\
          \bigjoinl_{j\in\n\setminus{\set i}} \preimage{\ff_{s}}{j}%
          \ar[r]%
          &\n\setminus{\set i}\ar[r,hookrightarrow]&%
          \numD{i-1}\join\preimage{\ff_{s'}}{j}\join\set{i+1,\dots,n}%
        \end{tikzcd}
      \end{equation*}
      The dashed arrow---which exists by the pushout property---exhibits
      $\preimage{\ff_{s'}}{i}$ as a retract of the singleton $\set{i}$, hence as
      a singleton itself.
    \item[~\ref{ax:gapsforpushout}]
      Fix $0<i\leq n$ and distinct elements $s,s'\in S$.
      Consider the commutative diagram
      \begin{equation*}
        \begin{tikzcd}
          \bigjoinl_{j\in\n}%
          \preimage{\ff_s}{j}\times%
          \preimage{\ff_{s'}}{j}%
          \ar[r]\ar[d]&%
          \bigjoinl_{j\in\n} \preimage{\ff_{s'}}{j}%
          \ar[d]\ar[rdd,dashed]\\%
          \bigjoinl_{j\in\n} \preimage{\ff_{s}}{j}%
          \ar[r]\ar[rrd,dashed]&%
          \n\ar[dr]\\
          &&\n
        \end{tikzcd}
      \end{equation*}
      where $\n\to\n$ is the (not order preserving) map
      that exchanges $i-1$ and $i$.
      By the pushout property of the solid square,
      at least one of the dashed composites must be not order preserving;
      this can only happen if least one of the maps $\ff_s$ and $\ff_{s'}$
      contains the subset $\set{i-1,i}\subseteq\n$ in its image.
      \qedhere
    \end{itemize}
  \end{enumerate}
\end{Prf}

\begin{Rem}
  \label{rem:pairwisecriterion}
  An $S$-claw $\exclaw=\tupleP{\ff_s}{s\in S}$ is backwards compatible if
  and only if for each pair of distinct elements $s,s'\in S$ the induced
  $\set{s,s'}$-subclaw is backwards compatible. Hence it follows
  from~\autoref{prop:coCartCechcubes}, that $\exclaw$ admits
  a \Cech{} cube in $\Deltaug$ if and only if each pair
  $\ff_s,\ff_{s'}$ (for distinct $s,s'\in S$) admits pullback in $\Deltaug$.
  Similarly, an $S$-claw admits a strongly biCartesian \Cech{} cube if and only
  if each two-pronged subclaw is compatible.
\end{Rem}

\subsection{Strongly biCartesian cubes in the cyclic category}

\newcommand{\cII}{I}

In this section, we characterize strongly biCartesian cubes in $\Lambda$.
To this end, we introduce the cyclic analog of a compatible claw.
Heuristically, this corresponds to adding the new \roughly{point} $(n,0)$
to the \roughly{manifold} $\n\in\Delta$.

\begin{Def}
  An $S$-claw $\exclaw\clawon\n$ in $\Delta$ is called
  \introduce{cyclically compatible} if
  the claw $\exclaw$ is compatible and all but at most one $\ff\in\exclaw$
  have the set $\set{0,n}\subseteq\n$ in their image.
\end{Def}

\begin{Rem}
  Let $\iota\colon\II''\to \II$ and $\alpha\colon \II''\to\II'$ be an inert map
  and an active map in $\Delta$, respectively.
  We can identify
  $\II=\II_0\join\II''\join\II_1$ and define $\n\coloneqq
  \II_0\join\II'\join\II_1$.
  It is easy to see that the $\numD 1$-claw
  $(\II'\hra \n, \Id\join\alpha\join\Id\colon\II\to\n)$
  is cyclically compatible and that $\II''$ is the associated pullback.
  By definition, the \buzzword{decomposition spaces} of
  \nameGCKT{}~\cite{GCKT2018a,GCKT2018b,GCKT2018c}
  are precisely those simplicial objects
  which send to Cartesian squares the biCartesian squares that arise this way.
\end{Rem}

\begin{Expl}
  \label{expl:decomp-claws}
  The $\numD 1$-claws
  \begin{equation}
    \begin{array}{*{4}c}
      0 & 1 & 2 & 3\\
      \hline
     \xxO&\xxX&\xxX&\xxX\\
     \xxX&\xxX&\xxO&\xxX
    \end{array}
    \intxt{and}
    \begin{array}{*{2}c}
      0 & 1 \\
      \hline
     \xxX&\xx2\\
     \xxO&\xxX
    \end{array}
  \end{equation}
  are cyclically compatible and arise as the pushouts
  of the inert map $\cofacemap{0}\colon\numD 1\to\numD 2$
  along the active maps
  $\cofacemap{1}\colon\numD 1\to\numD 2$ and
  $\codegmap{0}\colon\numD 1\to\0$, respectively.
  They encode the first upper $2$-Segal condition and an instance of unitality.
  The $\numD 1$-claw \eqref{eq:2} of \autoref{expl:1-Segal-claw}
  is not cyclically compatible because
  the \roughly{point} $(2,0)$ of the \roughly{manifold} $\numD 2$
  is not covered by any prong;
  the corresponding \Cech{} square \eqref{eq:3}
  is not coCartesian in the cyclic category.
\end{Expl}

The following is the main result of this section:

\begin{Prop}
  \label{prop:biCartinLambda}
  An $S$-claw $\exclaw\clawon\n$ in $\Delta$
  has a strongly biCartesian image in $\Lambda$
  if and only if it is cyclically compatible.
\end{Prop}

\begin{Cor}
  \label{cor:classesbiCartL}
  The following three classes of $S$-cubes in $\Lambda$ agree:
  \begin{itemize}
  \item
    strongly biCartesian $S$-cubes in $\Lambda$
  \item
    images of left active strongly biCartesian $S$-cubes in $\Delta$
  \item
    images of right active strongly biCartesian $S$-cubes in $\Delta$.
    \qedhere
  \end{itemize}
\end{Cor}

Before we can prove \autoref{prop:biCartinLambda}
and \autoref{cor:classesbiCartL} we need a couple of lemmas.

\begin{Lem}
  \label{lem:cyclic-compatible}
  Let $\exclaw=\tupleP{\ff_s\colon\II_s\to\n}{s \in S}$ be an $S$-claw in
  $\Delta$.
  If $\exclaw$ is compatible and either left active or right active then
  $\exclaw$ is cyclically compatible.
  Moreover, the following are equivalent:
  \begin{enumerate}
  \item
    \label{item:8}
    the claw $\exclaw$ is cyclically compatible;
  \item
    \label{item:9}
    for every $m\in\n$, the cyclic rotation
    $\cyclicrotation{+m}\exclaw\coloneqq \tupleP{
      \cyclicrotation{+m}{\ff_s}\colon
      \cyclicrotation{+m}{\II_s}\to\Dcyclicrotation{+m} n} {s\in S}$ of the
    claw $\exclaw$ is compatible;
  \item
    \label{item:10}
    there is an $m\in\n$ such that the cyclic rotation
    $\cyclicrotation{+m}\exclaw$ of the claw $\exclaw$ is left active and
    compatible;
  \item
    \label{item:11}
    there is an $m\in\n$ such that the cyclic rotation
    $\cyclicrotation{+m}\exclaw$ of the claw $\exclaw$ is right active and
    compatible.
    \qedhere
  \end{enumerate}
\end{Lem}

\begin{Prf}
  The first statement follows directly from the definitions.
  It is clear from the definition that the property of being
  cyclically compatible is preserved under cyclic rotation;
  hence we have the implications
  (\ref{item:8}$\implies$\ref{item:9}),
  (\ref{item:10}$\implies$\ref{item:8}) and
  (\ref{item:11}$\implies$\ref{item:8}).
  Given a compatible $S$-claw $\exclaw=\tupleP{\ff_s}{s\in S}$ on $\n$
  in $\Delta$, there is an element $m\in \n$ which is in the image of all the
  $\ff_s$.
  Then for any such $m$, the rotated claws
  $\cyclicrotation{-m}{\exclaw}$ and $\cyclicrotation{-m-1}\exclaw$
  are left active and right active, respectively.
  We thus obtain the implications
  (\ref{item:9}$\implies$\ref{item:10}) and
  (\ref{item:9}$\implies$\ref{item:11}).
\end{Prf}

\begin{Lem}
  \label{lem:CechcubesDL}
  Let $\excube\colon\powersetop S\to\Lambda$ be
  an $S$-cube in the cyclic category.
  The following are equivalent:
  \begin{enumerate}
  \item
    \label{item:1}
    the cube $\excube$ is strongly Cartesian;
  \item
    \label{item:2}
    there is a strongly Cartesian $S$-cube in $\Delta$
    which is mapped to $\excube$ under the canonical functor $\Delta\to\Lambda$;
  \item
    \label{item:4}
    every $S$-cube $\excube'$ in $\Delta$
    which maps to $\excube$ is strongly Cartesian.
    \qedhere
  \end{enumerate}
\end{Lem}

\begin{Prf}
  The implications
  \ref{item:2}$\implies$\ref{item:1}$\implies$\ref{item:4}
  follow from the general fact about slice categories that the projection
  $\Delta\cong\overcat\Lambda{\numL 0}\to\Lambda$
  preserves and reflects pullbacks.
  The implication
  \ref{item:4}$\implies$\ref{item:2}
  holds because the cube $\excube$ lifts to a cube in
  $\Delta\cong\overcat{\Lambda}{\numL0}$ by choosing any map
  $\excube(\emptyset)\to \numL 0$.
\end{Prf}

\begin{Lem}
  \label{lem:activepushoutsDL}
  Let
  \begin{equation}
    \label{eq:inlemsquare1}
    \begin{tikzcd}
      \II\cap\II'\rar{\ff\cap\II'}\ar[d,"\II\cap\ff'"']& \II'\dar{\ff'}\\
      \II\rar{\ff}&\n
    \end{tikzcd}
  \end{equation}
  be the left active strongly biCartesian \Cech{} square associated
  to a left active compatible claw $(\ff,\ff')\clawon\n$ in $\Delta$.
  Then the image in $\Lambda$
  of the square~\eqref{eq:inlemsquare1} is a pushout.
\end{Lem}

\begin{Prf}
  Consider a solid commutative diagram
  \begin{equation}
    \label{eq:inprfsquare2}
    \begin{tikzcd}
      \II\cap\II'\rar{\ff\cap\II'}\ar[d, "\II\cap\ff'"']&
      \II'\ar[rdd,"p'"]\ar[d,"\ff'"']\\
      \II\ar[rrd,"p"']\rar{\ff}&\numL n\ar[rd,dashed]\\
      &&\excycset
    \end{tikzcd}
  \end{equation}
  in $\Lambda$, where the top left square is the image
  of the square~\eqref{eq:inlemsquare1}.
  We need to show that there is a unique dashed morphism
  $p\colon\numL n\to\excycset$ of cyclic sets making the
  diagram~\eqref{eq:inprfsquare2} commute.

  \begin{itemize}
  \item
    First, we treat the case $\excycset=\numL0$. In this case the maps
    $p\colon \II\to\numL 0$, $p'\colon\II'\to\numL 0$ and
    $p''\colon\II\cap\II'\to\numL0$ correspond to cyclic rotations $\prec$ of
    the linear order on $\II$, $\II'$ and $\II''\coloneqq \II\cap\II'$,
    respectively; we have to show that there is a unique linear order $\prec$ on
    the cyclic set $\numL n$ such that both $\ff$ and $\ff'$ are order
    preserving with respect to $\prec$.
    Uniqueness is clear,
    because by compatibility of $(\ff,\ff')$
    each set $\set{i-1,i}$ (for $i\in\n$)
    is in the image of $\ff$ or of $\ff'$.

    To construct the linear order $\prec$ on $\n$,
    denote by $x$ and $x'$ the maximal elements in the linearly ordered
    sets $(\II,\prec)$ and $(\II',\prec)$, respectively, \ie the unique
    elements with $x+1\prec x$ and $x'+1\prec x'$. Without loss of generality,
    assume $i'\coloneqq\ff(x')\leq \ff(x)\eqqcolon i$. Define $\prec$ to be the
    unique linear order on the cyclic set $\numL n$ which has $i$ as its
    maximum. We need to show that $\ff$ and $\ff'$ preserve the orders $\prec$;
    for this it is enough to verify that
    $i<\ff(x+1)$ and $i<\ff'(x'+1)$
    (because $\ff(x)\leq i$ and $\ff'(x')\leq i$).
    
    Denote by $z''$, $z'$ and $z$ the $<$-minimal elements of $\II''$,
    $\II'$ and $\II$, respectively; they satisfy $(\ff\cap \II')(z'')=z'$,
    $(\II\cap\ff')(z'')= z$ and $\ff(z)= 0 =\ff'(z')$ because the
    square~\eqref{eq:inlemsquare1} was assumed to be left active.
    \begin{itemize}
    \item Assume that $i=\ff(x)=\ff(x+1)$.
      Then by backwards compatibility of $(\ff,\ff')$
      we must have a unique $y'\in\II'$ with $\ff'(y')=i$.
      By the explicit formula for \Cech{} cubes
      we deduce that the order preserving map
      (with respect to both $\prec$ and $<$)
      $\II\cap \ff'\colon\II''\to\II$
      restricts to a bijection
      $\II''\cap\set{i}\xra{\cong}\II\cap \set{i}$
      which is therefore an isomorphism
      (with respect to $\prec$ and $<$).
      Denote by $\ol x, \ol{x+1}\in\II''$ the (unique) preimages under
      $\II\cap\ff'$ of $x$ and $x+1$, respectively;
      they satisfy $\ol x +1= \ol{x+1}\prec\ol{x}$ by the isomorphism property,
      which means they are the maximal and minimal element of the linearly
      ordered set $(\II'',\prec)$, respectively.
      Since both $\ol x$ and $\ol{x+1}$ are mapped to $y'$
      by $\ff\cap\II'$ we deduce that
      $\ff\cap\II'\colon \II''\to\II'$ is constant.
      This can only happen if
      $\ff$ was already constant and $\ff'$ was an equivalence.
      Hence the square~\eqref{eq:inlemsquare1} is degenerate
      and therefore trivially a pushout in $\Lambda$.
    \item The case $i'=\ff'(x')=\ff(x'+1)$ is analogous.
    \end{itemize}
    We may therefore assume that $x$ and $x'$ are the maximal elements
    (with respect to both $<$ and $\prec$)
    of their corresponding preimages
    $\preimage{\ff} i$ and $\preimage{\ff'}{i'}$.
    It follows directly that
    $\ff(x+1)>i$ and $\ff'(x'+1)>i'$;
    it remains to show $\ff'(x'+1)>i$ and we may assume that $i'<i$.
    Next, we show that there is no $j\in\n$ with $i'<j\leq i$
    which is in the image of
    $\ff''\coloneqq\ff\cap\ff'\colon\II''\to\n$:
    \begin{itemize}
    \item Otherwise, choose $w''\in\II''$ with $\ff''(w'')=j$.
      Set $w'\coloneqq (\ff\cap\II')(w'')\in\II'$
      and $w\coloneqq(\II\cap\ff')(w'')\in\II$.
      We have $z<w$ and $z'\leq x'<w'$ by construction
      and $w\leq x$ because $x$ is maximal for $<$
      in the preimage $\preimage\ff{i}$.
      Hence we have
      (after cyclic rotation and using that $x$ and $x'$ are $\prec$-maximal)
      $z\prec w\preceq x$ and $w'\prec z'\preceq x'$,
      which implies $z''\prec w''$ and $w''\prec z''$, respectively.
      Contradiction.
    \end{itemize}
    Since $i$ is in the image of $\ff$ (by definition) and
    each $j$ with $i'<j\leq i$ is not in the image of $\ff''$,
    it follows from the compatibility of $(\ff,\ff')$
    that each such $j$ is not in the image of $\ff'$.
    Since we already know $\ff'(x'+1)>i'$
    we obtain $\ff'(x'+1)>i$, as desired;
    this concludes the case $\excycset=\numL 0$.
  \item
    We prove the case of a general $\excycset$.
    To see the existence of the dashed map
    in the diagram~\eqref{eq:inprfsquare2},
    choose any map $\excycset\to\numL0$.
    By the case $\excycset=\numL 0$ which we have just shown,
    we can fill the dotted morphism
    $\numL n\to\numL 0$
    of cyclic sets in the following commutative diagram
    \begin{equation}
      \label{eq:inprfsquare3}
      \begin{tikzcd}
        \II\cap\II'\rar{\ff\cap\II'}\ar[d, "\II\cap\ff'"']&
        \II'\ar[rdd,"p'"]\ar[d,"\ff"']\\
        \II\ar[rrd,"p"']\rar{\ff}&\numL n\ar[rd,dashed]\ar[rrd,dotted]\\
        &&\excycset\ar[r]&\numL 0
      \end{tikzcd}
    \end{equation}
    Thus we have constructed a diagram in the overcategory
    $\overcat\Lambda{\numL 0}$.
    Under the canonical identification
    $\Delta\cong\overcat\Lambda{\numL 0}$,
    the top left square of the
    diagram~\eqref{eq:inprfsquare3}
    gets identified with a cyclic rotation of
    the original diagram~\eqref{eq:inlemsquare1}.
    Since any cyclic rotation of a left active compatible claw is compatible,
    we deduce from \autoref{cor:biCart=comp}
    that the corresponding \Cech{} square is a pushout in
    $\Delta\cong\overcat\Lambda{\numL 0}$.
    We conclude by the pushout property that
    the desired dashed map $\numL n\to\excycset$
    in~\eqref{eq:inprfsquare3} and
    \latin{a fortiori} in~\eqref{eq:inprfsquare2} exists.

    To prove uniqueness,
    recall that the square~\eqref{eq:inlemsquare1} is a pushout
    on the level of underlying sets,
    so that the dashed map is unique as a function of underlying sets.
    If $\numL n\to\excycset$ is constant then
    it factors uniquely as $\numL n\to\numL 0\to\excycset$,
    hence is unique by the case $\excycset=\numL 0$.
    If $\numL n\to\excycset$ is not constant then
    it is uniquely determined by its underlying function of sets.
    \qedhere
  \end{itemize}
\end{Prf}

\begin{Prf}[of \autoref{prop:biCartinLambda}]
  If $\exclaw$ is cyclically compatible then by \autoref{lem:cyclic-compatible}
  there is a cyclic rotation $\cyclicrotation{-m}\exclaw$ of
  $\exclaw$ which is left active and compatible.
  Since $\exclaw$ and $\cyclicrotation{-m}{\exclaw}$ have isomorphic images in
  $\Lambda$, it is enough to show that the latter image is strongly biCartesian.
  Since the \Cech{} cube $\Cechcube{\cyclicrotation{-m}\exclaw}$
  is left active and strongly biCartesian,
  it follows from \autoref{lem:CechcubesDL} and \autoref{lem:activepushoutsDL}
  (applied to each $2$-dimensional face of the cube)
  that its image in $\Lambda$ is still strongly biCartesian.
   
  Conversely,
  let $\excube$ be a strongly biCartesian cube in $\Lambda$
  extending $\exclaw$.
  Then every choice of $m\in\n$ yields a structure map
  $\Dcyclicrotation {+m}n\colon\excube(\emptyset) =\numL n\to \numL 0$
  which gives rise to a cube $\excube_m$
  in $\overcat\Lambda{\numL 0}\cong \Delta$
  that maps to $\excube$ and extends the claw $\cyclicrotation{+m}\exclaw$.
  Since the slice projection $\Delta\to\Lambda$ reflects
  pullbacks and pushouts,
  we deduce that each of these cubes $\excube_m$ is strongly
  biCartesian. Hence by \autoref{cor:biCart=comp} the corresponding claw
  $\cyclicrotation{+m}\exclaw$ is compatible. We conclude
  by \autoref{lem:cyclic-compatible} that the original claw $\exclaw$ is
  cyclically compatible.
\end{Prf}

\begin{Prf}[of \autoref{cor:classesbiCartL}]
  Recall from \autoref{cor:biCart=comp} that strongly biCartesian $S$-cubes in
  $\Delta$ are precisely the \Cech{} cubes of compatible $S$-claws.
  Hence \autoref{cor:classesbiCartL} follows directly from
  \autoref{prop:biCartinLambda} and \autoref{lem:cyclic-compatible}.
\end{Prf}

\subsection{Primitive decomposition of biCartesian cubes}

In this section we show how a strongly biCartesian cube in $\Delta$
can be decomposed into simpler building blocks.

\begin{Def}
  A map $\ff\colon\II\to\n$ in $\Delta$ is called \introduce{primitive}
  if there is exactly one $i\in\n$
  such that $\preimage\ff i$ is not a singleton;
  the map $\ff$ is called \introduce{preprimitive}
  if it is primitive or an isomorphism.
  A candidate covering $\exclaw$ in $\Deltaug$
  (and the corresponding \Cech{} cube $\Cechcube\exclaw$)
  is called \introduce{(pre)primitive}
  if the claw $\exclaw$ consists only of (pre)primitive maps.
\end{Def}

\newcommand{\canfactterm}[2]{#1_{#2}}
\newcommand{\canfactcomp}[2]{#1_{#2}}
\newcommand{\canfactmap}[2]{\ol {#1}_{#2}}
\newcommand{\iin}{i}

\begin{Cstr}
  \label{cstr:canonicalfactorization}
  Let $\ff\colon \II\to\n$ be a map in $\Delta$.
  For each $\iin\in\set{-1,0,\dots,n}$, we define objects
  \begin{equation*}
    \canfactterm \II \iin\coloneqq
    \invimage\ff{\numD{\iin}}
    \join \numDminus{n}{\iin}
  \end{equation*}
  in $\Delta$.
  Then $\ff$ admits a canonical factorization
  \begin{equation}
    \label{eq:canfact}
    f\colon \II=\canfactterm\II n\xra{\canfactmap\ff n}\
    \dots\xra{\canfactmap\ff {\iin+1}}
    \canfactterm \II \iin\xra{\canfactmap \ff \iin}
    \dots \xra{\canfactmap\ff 1}
    \canfactterm \II{0}\xra{\canfactmap \ff 0}\canfactterm \II{-1}=\n
  \end{equation}
  where each map
  $\canfactmap \ff \iin\colon\canfactterm \II i\to \canfactterm \II{\iin-1}$
  is given as
  \begin{equation*}
    \canfactmap\ff \iin\coloneqq \Id_{\invimage\ff{\numD{\iin-1}}}%
    \join\left(
      \ff\cap\set{\iin}\colon\preimage \ff \iin\to\set{\iin}\right)%
    \join\Id_{\numDminus{n}{\iin}}.%
  \end{equation*}
  Observe that each map $\canfactmap{\ff}\iin$ is preprimitive.
\end{Cstr}

\newcommand{\Cechcubeindirection}[2]{\check{\mathrm C}^{#2}{#1}}
\begin{Lem}
  \label{lem:liftfactorization}
  Let ($\ff\colon \II\to\n, \ff'\colon\II'\to\n)$
  be backwards compatible and factorize $\ff$ as in
  \autoref{cstr:canonicalfactorization}.
  \begin{enumerate}
  \item
    For each $\iin\in\n$, the composition
    $\canfactterm{\II}{\iin}\to\n$
    in \eqref{eq:canfact}
    is backwards compatible with $\ff'$ so that
    by~\autoref{prop:coCartCechcubes} we can form the pullbacks
    \begin{equation}
      \begin{tikzcd}
        \II\cap\II'\dar\rar\ar[d,"\II\cap\ff'"'] &%
        \canfactterm{\II}{n-1}\cap\II'\rar\ar[d,"\II_{n-1}\cap\ff'"']&%
        \dots\rar&%
        \canfactterm{\II}{1}\cap\II'\rar\ar[d,"\II_1\cap\ff'"]&%
        \canfactterm{\II}{0}\cap\II'\rar\ar[d,"\II_0\cap\ff'"]&%
        \II'\dar{\ff'}\\%
        \II\rar{\canfactmap \ff n}\ar[rrrrr,"\ff"',bend right=10]&%
        \canfactterm \II {n-1}\ar[r,"\canfactmap \ff {n-1}"]&%
        \dots\ar[r,"\canfactmap\ff 2"]&%
        \canfactterm \II 1\rar{\canfactmap \ff 1}&%
        \canfactterm \II{0}\rar{\canfactmap \ff 0}&%
        \n
      \end{tikzcd}
    \end{equation}
    which factorize the \Cech{} square of $\ff$ and $\ff'$
    into smaller \Cech{} squares.
  \item
    The original claw $(\ff,\ff')$ is compatible if and only if
    the claw
    $(\canfactmap{\ff}{\iin},\canfactterm{\II}{\iin-1}\cap\ff')
    \clawon\canfactterm{\II}{\iin-1}$
    is compatible
    for each $\iin\in \n$.
  \item
    The original claw $(\ff,\ff')$ is cyclically compatible if and only if
    the claw
    $(\canfactmap{\ff}{\iin},\canfactterm{\II}{\iin-1}\cap\ff')
    \clawon\canfactterm{\II}{\iin-1}$ is cyclically compatible
    for each $\iin\in \n$.
    \qedhere
  \end{enumerate}
\end{Lem}

\begin{Prf}
  Follows by direct inspection of the explicit constructions.
\end{Prf}

\begin{Lem}
  \label{cor:decompintoprimitives}
  \begin{enumerate}
  \item
    \label{item:3}
    Every strongly biCartesian cube $\excube$ in $\Delta$ can
    be decomposed into a pasting of preprimitive strongly biCartesian cubes.
    If $\excube$ was left active then each of these cubes can
    be chosen to be left active.
    If $\excube$ was right active then each of these cubes can be chosen to be
    right active.
  \item
    \label{item:5}
    Every cube in $\excube$ in $\Delta$
    which becomes strongly biCartesian in $\Lambda$ can be
    decomposed into a pasting of preprimitive strongly biCartesian cubes,
    each of which is left active or right active.
  \item
    \label{item:6}
    If the original cube $\excube$ in
    \ref{item:3} or \ref{item:5} is nondegenerate
    then the pastings can be chosen to consist of primitive cubes.
    \qedhere
  \end{enumerate}
\end{Lem}

\begin{Prf}
  By \autoref{cor:biCart=comp}, each strongly biCartesian cube in $\Delta$
  is the \Cech{} cube $\Cechcube\exclaw$
  of some compatible $S$-claw $\exclaw=\tupleP{\ff_s}{s\in S}$.
  By \autoref{prop:biCartinLambda}, each cube in $\Delta$ which becomes
  strongly biCartesian in $\Lambda$ is of this form $\Cechcube\exclaw$ where
  $\exclaw$ is cyclically compatible.
  For each $s\in S$,
  consider the factorization of $\ff_s$ into preprimitive maps
  from \autoref{cstr:canonicalfactorization}.
  By a repeated application of \autoref{lem:liftfactorization},
  we can decompose the cube $\Cechcube\exclaw$ into a pasting of \Cech{} cubes
  of compatible claws which are cyclically compatible if $\exclaw$ was.
  Parts~\ref{item:3} and~\ref{item:5} of \autoref{cor:decompintoprimitives}
  now follow by applying
  \autoref{cor:biCart=comp}, \autoref{prop:biCartinLambda}
  and by the observing that preprimitive cyclically compatible claws
  are automatically either left active or right active.
  Part~\ref{item:6} follows with the same procedure by
  dropping all identities appearing in the factorizations
  produced by \autoref{cstr:canonicalfactorization}.
\end{Prf}

\section{Precovers and intersection cubes}
\label{sec:prec-inters-cubes}

Let $\exclaw\clawon\n$ be a $S$-claw on $\n$ in $\Delta$.
If all of the maps in the claw $\exclaw$ are injective
then we call $\exclaw$ an \introduce{($S$-)precover} on $\n$.
Since precovers are trivially backwards compatible,
\autoref{prop:coCartCechcubes} guarantees the existence
of the \Cech{} cube $\Cechcube\exclaw$;
we call it the \introduce{intersection cube} of $\exclaw$.
If we view the injective maps $\exclaw\ni\ff_s\colon\II_s\hra\n$ as subsets
$\II_s\subseteq\n$ of $\n$
then the intersection cube of $\exclaw$ is given explicitly by the intersections
\begin{equation}
  \label{eq:explicitintercube}
  T\lmapsto \bigcap_{t\in T}\II_t,
\end{equation}
(where the empty intersection is $\n$ by convention);
thus the terminology \qquote{intersection cube} is justified.

\subsection{Membrane spaces and refinements}
\label{sec:membrane-spaces}

By right Kan extension along the Yoneda embedding
$\Delta\hra\Fun(\Dop,\Set)$,
we can extend any simplicial object $\xdc$ to a functor
\begin{equation*}
  \localX\colon\Fun(\Dop,\Set)^\op\lra\localC,
\end{equation*}
which we still denote by $\localX$.
Given any simplicial set $\exsset$,
we can calculate the value of $\localX$ at $\exsset$%
---which \twonames{Dyckerhoff}{Kapranov} call the
\introduce{object of $\exsset$-membranes in $\localX$}%
---by the pointwise formula for Kan extensions:
\begin{equation*}
  \Xonsset\exsset\simeq
  \lim\left(\left(\overcat\Delta\exsset\right)^\op
    \to\Dop\xra{\localX}\localC\right)
\end{equation*}
The inclusion $\Delta\hra\Fun(\Dop,\Set)$ factors as
$\Delta\hra\Deltaug\hra\Fun(\Dop,\Set)$,
where the second map sends the initial object $\emptyset$
to the initial presheaf.
We can therefore evaluate any simplicial objet $\xdc$ at $\emptyset$ and the
value will be a terminal object in $\localC$.

Given a candidate covering
$\exclaw=\tupleP{\ff_s\colon\II_s\to\n}{s\in S}$ in $\Delta$,
we obtain a simplicial set $\ssetofprecover{\exclaw}$ as the colimit
\begin{equation*}
  \ssetofprecover{\exclaw}\coloneqq
  \colim
  \left(
    \powersetneop{S}
    \xra{\Cechcube\exclaw}
    \Delta
    \lhra
    \Fun(\Dop,\Set)
  \right)
\end{equation*}
which comes equipped with a canonical map
$\ssetofprecover{\exclaw}\to \simplex n$.
It is easy to see that if $\exclaw$ is a precover
(\ie if all maps $\ff_s$ are injective)
then
$\ssetofprecover\exclaw\subseteq \simplex n$
can be identified with the simplicial subset
$\ssetofprecover\exclaw\coloneqq\bigcup_{\II_s\in S}\simplex {\II_s}$
of the $n$-simplex.
We say that a precover $\exclaw'\clawon\n$
is a \introduce{refinement} of $\exclaw\clawon\n$---%
written $\exclaw'\refines\exclaw$---%
if and only if $\ssetofprecover{\exclaw'}$
is a simplicial subset of $\ssetofprecover{\exclaw}$;
explicitly, this means that
for every $\II'\in\exclaw'$
there is at least one $\II\in\exclaw$
such that $\II'\subseteq \II$ (as subobjects of $\n$).
We say the refinement $\exclaw'\refines\exclaw$ is
\introduce{degenerate}
if $\ssetofprecover{\exclaw'}=\ssetofprecover{\exclaw}$.
For each $\n\in\Delta$ the assignment
$\exclaw\mapsto\ssetofprecover{\exclaw}$
describes an equivalence of categories between
the category (which is just a preorder) of precovers and refinements on $\n$
and the full subcategory of the overcategory
$\overcat{\Fun(\Dop,\Set)}{\simplex n}$
spanned by the simplicial subsets of $\simplex n$.
An explicit inverse is given by identifying each simplicial subset
$\exsset\subseteq \simplex n$
with the precover given by the maximal simplices of $\exsset$.
We will implicitly use this identification and write
\begin{equation*}
  \reducedclaw{\exclaw}\coloneqq
  \tupleP{\II}{\simplex\II\hra\ssetofprecover{\exclaw}\text{ maximal}}
  \clawon\n
\end{equation*}
for the precover obtained from a precover $\exclaw$ by
\roughly{removing redundant subsets}.

\begin{Rem}
  \label{rem:identification-cubes}
  For every precover $\exclaw$,
  the restriction
  $\restr{\Cechcube\exclaw}{\powersetneop S}\colon
  \powersetneop{S}\to \overcat\Deltaug{\ssetofprecover\exclaw}$
  of the \Cech{} cube of $\exclaw$
  has a left adjoint given by
  \begin{equation*}
    (\m, \alpha\colon\simplex m \to \ssetofprecover\exclaw)
    \lmapsto
    \setP{s\in S}
    {\alpha (\simplex m)\subseteq\simplex{\II_s}}
  \end{equation*}
  Since right adjoints are homotopy initial\footnote{
    Here we use the terminology of Dugger~\cite{Dugger}:
    He calls homotopy terminal what Joyal and Lurie
    would call cofinal;
    homotopy initial is then the dual notion.
  },
  the canonical map
  \begin{equation*}
    \Xonsset{\ssetofprecover\exclaw}
    \simeq
    \lim\restr{\localX}{\left(
        \overcat\Deltaug{\ssetofprecover\exclaw}\right)^\op}
    \xra{\simeq}
    \lim\restr{\localX\circ \Cechcube\exclaw}{\powersetne S}
  \end{equation*}
  is an equivalence.
  In particular,
  $\localX$ satisfies descent with respect to $\exclaw$
  if and only if
  $\localX$ sends the inclusion
  $\ssetofprecover{\exclaw}\hra\simplex n$ to an equivalence.
\end{Rem}

\begin{Def}
  We say that a refinement $\exclaw'\refines\exclaw$
  of precovers $\n$ is $\localX$-local if the induced morphism
  $\ssetofprecover{\exclaw'}\to\ssetofprecover{\exclaw}$
  of simplicial sets is sent by $\localX$ to an equivalence in $\localC$.
\end{Def}

The following lemma (which is essentially Corollary 3.16 in~\cite{DJW2018})
is the main tool to compare to one another
descent conditions with respect to various precovers.

\begin{Lem}
  \label{lem:one-step-refinement}
  Let $\exclaw\clawon\n$ be a precover in $\Delta$ and
  $\II\subset \n$ a subset.
  Assume that the restricted precover
  \begin{equation*}
    \exclaw\cap\II\coloneqq\tupleP{\II'\cap\II}{\II'\in\exclaw}\clawon\II
  \end{equation*}
  on $\II$ is $\localX$-local.
  Then the refinement
  $\exclaw\refines \reducedclaw{\exclaw\cup\set{\II}}$ is $\localX$-local.
  In particular,
  the original precover $\exclaw$ is $\localX$-local
  if and only if
  the extended precover $\reducedclaw{\exclaw\cup{\II}}$ is $\localX$-local.
\end{Lem}

\begin{Prf}
  The refinement $\exclaw\refines\reducedclaw{\exclaw\cup\set{\II}}$
  can be written as the composition of refinements
  \begin{equation}
    \label{eq:inprftisntrise}
    \exclaw\refines\exclaw\cup\set{\II}
    \refines\reducedclaw{\exclaw\cup\set{\II}}.
  \end{equation}
  The first refinement in the composition~\eqref{eq:inprftisntrise}
  is $\localX$-local by \autoref{lem:changeofclaw}
  (due to the assumption of the lemma
  and using the identification of \autoref{rem:identification-cubes});
  the second refinement is degenerate, hence always local.
  The claim follows.
\end{Prf}

\subsection{Polynomial simplicial objects}
\label{sec:polyn-simpl-objects}

Recalling the analogy to manifold calculus
described in \autoref{sec:manifold-calc-Delta},
we observe that compatible precovers can be identified precisely
with the \roughly{open covers} of the form~\eqref{eq:1}.
Indeed, an $S$-precover $\exclaw$ on $\n\in \Delta$ is compatible
if and only if every \roughly{point} $(x-1,x)$ of the \roughly{manifold} $\n$
is contained in all but at most one of the elements of $\exclaw$,
which we think of as \roughly{open subsets} of $\n$;
in other words, $\exclaw$ consists precisely of \roughly{open subsets}
with \roughly{pairwise disjoint closed complements}.
The analogy thus motivates the following definition:

\begin{Def}
  \label{def:polynomial-f}
  We call a functor $\Dop\to\localC$
  \introduce{polynomial of degree $\leq\card{S}$}
  (or \introduce{$S$-polynomial}, for short)
  if $\localX$ satisfies descent with respect
  to all compatible $S$-covers in $\Delta$.
\end{Def}

\begin{Expl}
  We depict, for $\kk=1,2,3$,
  the unique nondegenerate compatible $\numD\kk$-cover on $\numD{2\kk}$:
  \begin{equation}
    \label{eq:123comp-claw}
    \begin{array}{*{3}c}
      0 & 1 & 2\\
      \hline
      \xxO&\xxX&\xxX\\
      \xxX&\xxX&\xxO
    \end{array}
    \intxt{}
    \begin{array}{*{5}c}
      0 & 1 & 2 & 3 & 4\\
      \hline
      \xxO&\xxX&\xxX&\xxX&\xxX\\
      \xxX&\xxX&\xxO&\xxX&\xxX\\
      \xxX&\xxX&\xxX&\xxX&\xxO
    \end{array}
    \intxt{}
    \begin{array}{*{7}c}
      0 & 1 & 2 & 3 & 4 & 5 & 6\\
      \hline
      \xxO&\xxX&\xxX&\xxX&\xxX&\xxX&\xxX\\
      \xxX&\xxX&\xxO&\xxX&\xxX&\xxX&\xxX\\
      \xxX&\xxX&\xxX&\xxX&\xxO&\xxX&\xxX\\
      \xxX&\xxX&\xxX&\xxX&\xxX&\xxX&\xxO
    \end{array}
  \end{equation}
  Note that for $n<2\kk$,
  there are no nondegenerate compatible $\numD\kk$-covers on $\n$.
\end{Expl}

The number of compatible $S$-covers on $\n\in\Delta$
grows quite rapidly in $n$.
Thus \latin{a priori} to determine that a simplicial object is $S$-polynomial,
there is an increasing number of conditions to check in each dimension.
We show now that it suffices to check \emph{any one}
non-trivial condition in each dimension.

\begin{Prop}
  \label{lem:onebiCartforall}
  Let $\xdc$ be a simplicial object
  in some $\infty$-category with finite limits.
  Assume that for each $n\geq 2\kk$ there exists
  a nondegenerate compatible $\numD\kk$-cover $\exclaw\clawon\n$ in $\Delta$
  which is $\localX$-local.
  Then all compatible $\numD\kk$-covers are $\localX$-local.
\end{Prop}

\begin{Prf}
  Assume the assumption of \autoref{lem:onebiCartforall}.
  Recall that degenerate covers are automatically local.
  Hence there is nothing to show for $n<2\kk$
  because in this case there are
  no nondegenerate compatible $\numD\kk$-covers on $\n$.
  We prove by induction on $n\geq 2\kk$ that
  all nondegenerate compatible $\numD\kk$-covers are $\localX$-local.
  The inductions start is the case $n= 2\kk$,
  which is trivial because
  there is a unique nondegenerate compatible $\numD\kk$-cover
  on $\numD{2\kk}$.
  For the induction step consider the following directed graph:
  \begin{itemize}
  \item
    Vertices are nondegenerate compatible $\numD\kk$-covers on $\n$.
  \item
    Let $\exclaw$ be a nondegenerate compatible $\numD\kk$-cover
    and let $\II\in\exclaw$ and $x\in\n\setminus{\II}$ such that
    $\II'\coloneqq\II\cup\set{x}\neq \n$.
    Then the cover
    $\exclaw'\coloneqq\reducedclaw{\exclaw\cup{\II'}}$
    is easily seen to be again
    $\numD\kk$-pronged, compatible and nondegenerate.
    We add the refinement
    \begin{equation*}
      \exclaw \refines \reducedclaw{\exclaw\cup{\II'}}
    \end{equation*}
    to the graph as an arrow $\exclaw\to\exclaw'$.
    Observe that in the language of \autoref{rem:visualization},
    the cover $\exclaw'$ arises from the cover $\exclaw$
    by choosing a row with at least two $\xxO$'s
    and replacing one of them by $\xxX$.
  \end{itemize}
  With the notation above it is easy to see that
  the restricted $\numD\kk$-cover
  $\exclaw\cap\II'\clawon\II'$ is still compatible,
  hence $\localX$-local by the induction hypothesis
  (since $\II'\subsetneq \n$).
  It follows from~\autoref{lem:one-step-refinement} that every arrow in
  the graph corresponds to an $\localX$-local refinement.
  The proof of \autoref{lem:onebiCartforall}
  is concluded by the easy combinatorial observation that
  the graph is connected as an undirected graph,
  \ie one can connect every pair of nondegenerate compatible $\numD\kk$-covers
  by a zigzag of $\localX$-local refinements as above.
\end{Prf}

\begin{Rem}
  The directed graph constructed in the proof of \autoref{lem:onebiCartforall}
  is just the Hasse diagram of the poset
  of nondegenerate compatible $\numD\kk$-covers under refinement.
  Our proof therefore shows that if there is an $n\geq 2k$
  such that $\localX$ satisfies descent with respect to all compatible
  $\numD\kk$-covers in $\Deltall{n}$
  then all refinements between nondegenerate compatible $\numD\kk$-covers on
  $\n$ are $\localX$-local.
\end{Rem}

\section{Weakly excisive and weakly $\Lambda$-excisive simplicial objects}
\label{sec:excisive-objects}

Fix an $\infty$-category $\localC$ with finite limits.
Recall from \autoref{sec:weak-excision} that a simplicial object $\xdc$ is
\begin{itemize}
\item
  weakly $S$-excisive if it sends strongly biCartesian $S$-cubes in $\Delta$
  to Cartesian cubes in $\localC$.
\item
  weakly $S$-$\Lambda$-excisive if it sends to Cartesian cubes in $\localC$
  those $S$-cubes in $\Delta$ which become strongly biCartesian in $\Lambda$
  after applying the canonical functor $\Delta\to\Lambda$.
\end{itemize}

\begin{Rem}
  \label{rem:cyclic-excisive}
  It follows from \autoref{rem:restr-relative-exc}
  that a cyclic object $\Lambda^\op\to\localC$ is weakly $S$-excisive
  if and only if its restriction to $\Delta$ is weakly $S$-$\Lambda$-excisive.
\end{Rem}

We can refine the notion of weak $\Lambda$-excision as follows:

\begin{Def}
  A simplicial object $\xdc$ in $\localC$ is called
  \begin{itemize}
  \item
    \introduce{lower weakly $S$-$\Lambda$-excisive}
    if $\localX$ sends
    every left active strongly biCartesian $S$-cube in $\Delta$
    to a Cartesian cube in $\localC$;
  \item
    \introduce{upper weakly $S$-$\Lambda$-excisive}
    if $\localX$ sends
    every right active strongly biCartesian $S$-cube in $\Delta$
    to a Cartesian cube in $\localC$.
    \qedhere
  \end{itemize}
\end{Def}

The terminology is justified by the following easy lemma.

\begin{Lem}
  A simplicial object is weakly $S$-$\Lambda$-excisive if and only if it
  is both lower weakly $S$-$\Lambda$-excisive and upper weakly $S$-$\Lambda$-excisive.
\end{Lem}

\begin{Prf}
  By \autoref{cor:decompintoprimitives},
  every $S$-cube in $\Delta$ with strongly biCartesian image in $\Lambda$
  can be decomposed into a pasting of strongly biCartesian cubes
  each of which is left active or right active;
  thus we have \qquote{if}.
  The converse \qquote{only if} follows from the fact
  (\autoref{cor:classesbiCartL})
  that every strongly biCartesian in $\Delta$
  which is left active or right active
  has a strongly biCartesian image in $\Lambda$.
\end{Prf}

\subsection{Weakly Excisive = polynomial}

As explained in \autoref{sec:polyn-simpl-objects},
a polynomial functor of degree $\geq\kk$
is a simplicial object $\Dop\to\localC$ which sends
all strongly biCartesian intersection $\numD\kk$-cubes
to Cartesian cubes in $\localC$.
\latin{A priori}, this does not agree with weak $\kk$-excision,
because it only takes into account
strongly biCartesian cubes which consist of \emph{injective} maps.
The next theorem states that this discrepancy is illusory
both for weak ($\Delta$-)\-excision and for
(lower and/or upper) weak $\Lambda$-excision.

\begin{Thm}
  \label{thm:unitality-of-cubes}
  Let $\localC$ be an $\infty$-category with all finite limits.
  A simplicial object $\xdc$ is
  \begin{enumerate}[label=(\alph*)]
  \item
    \label{item:excisiveinter}
    weakly $S$-excisive if and only if it sends primitive strongly biCartesian
    intersection $S$-cubes in $\Delta$ to Cartesian cubes in $\localC$;
  \item
    \label{item:Lexcisiveinter}
    lower weakly $S$-$\Lambda$-excisive if and only if it sends primitive strongly
    biCartesian left active intersection $S$-cubes in $\Delta$ to Cartesian
    cubes in $\localC$;
  \item
    \label{item:Rexcisiveinter}
    upper weakly $S$-$\Lambda$-excisive if and only if it sends primitive strongly
    biCartesian right active intersection $S$-cubes in $\Delta$ to Cartesian
    cubes in $\localC$.
    \qedhere
  \end{enumerate}
\end{Thm}

Before we prove \autoref{thm:unitality-of-cubes},
we deduce the following criterion for detecting
weak $\Lambda$-excision of a simplical object in terms of
weak ($\Delta$-)excision of its path objects.

\begin{Cor}[Path space criterion]
  \label{cor:path-space-criterion-cubes}
  A simplicial object $\xdc$ in an $\infty$-category with all finite limits is
  \begin{itemize}
  \item
    lower weakly $S$-$\Lambda$-excisive if and only if the left path object
    $\lpath\localX\coloneqq\localX\circ(\0\join\blank)$ is weakly $S$-excisive;
  \item
    upper weakly $S$-$\Lambda$-excisive if and only if the right path object
    $\rpath\localX\coloneqq\localX\circ(\blank\join\0)$ is weakly $S$-excisive.
    \qedhere
  \end{itemize}
\end{Cor}

\begin{Prf}
  Observe that composition with the functor
  $\0\join\blank\colon\Delta\to\Delta$ identifies
  compatible $S$-covers in $\Delta$
  with left active compatible $S$-covers in $\Delta$;
  hence by \autoref{cor:biCart=comp} it identifies
  strongly biCartesian intersection $S$-cubes in $\Delta$
  with left active strongly biCartesian intersection $S$-cubes $\Delta$.
  The first statement of \autoref{cor:path-space-criterion-cubes}
  now follows directly from \autoref{thm:unitality-of-cubes};
  the proof of the second statement is analogous.
\end{Prf}

\begin{Rem}
  The proof of \autoref{cor:path-space-criterion-cubes}
  makes crucial use of \autoref{thm:unitality-of-cubes}
  because in general a left active diagram in $\Delta$
  need not factor through the functor $\0\join\blank\colon\Delta\to\Delta$.
  It is the fact that we can reduce to diagrams of \emph{injective} maps
  that makes this argument work.
\end{Rem}

To prove \autoref{thm:unitality-of-cubes}
we isolate the following key lemma which we prove separately below.
Recall that, for each $m\geq0$, we denote the unique active maps
$\numD 1\to\numD m$ in $\Delta$ by $\canactivemap{m}$.

\begin{Lem}[Key lemma]
  \label{lem:activetozero}
  Let $p\colon\localC\to \calB$ be a Cartesian fibration of
  $\infty$-categories. Let $\xdc$ be a simplicial object. Assume that, for all
  $m\geq 1$, the edge $\localX(\canactivemap m)$ of $\localC$ is $p$-Cartesian.
  Then the edge $\localX(\alpha)$ is also $p$-Cartesian for every active
  morphism $\alpha$ in $\Delta$.
\end{Lem}

\begin{Prf}[of \autoref{thm:unitality-of-cubes}]
  We will prove part~\ref{item:excisiveinter};
  the proof for \ref{item:Lexcisiveinter} or \ref{item:Rexcisiveinter}
  is the same, word by word, by only considering cubes
  which are left or right active, respectively.
  The direction \qquote{only if} is trivial.
  
  To prove \qquote{if}
  let $\xdc$ be a simplicial object which sends
  primitive strongly biCartesian intersection $S$-cubes in $\Delta$
  to Cartesian cubes in $\localC$.
  Assume that there is a counterexample to \autoref{thm:unitality-of-cubes},
  \ie a compatible $S$-claw $\exclaw=\tupleP{\ff_s}{s\in S}$ on $\n\in\Delta$
  such that the corresponding \Cech{} cube $\Cechcube\exclaw$
  is not sent by $\localX$ to a Cartesian cube in $\localC$.
  By \autoref{cor:decompintoprimitives} we
  may choose $\exclaw$ to be preprimitive.
  We may assume that $\exclaw$ is primitive because otherwise it would be
  degenerate; and degenerate cubes are always sent to Cartesian cubes.
  By induction we may additionally assume that the number
  \newcommand\localdegeneracycount[1]{d{#1}}
  \begin{equation}
    \localdegeneracycount{\exclaw}\coloneqq
    \card{\setP{s\in S}{\ff_s\text{ is not injective}}}
  \end{equation}
  is minimal amongst all counterexamples.
  The number $\localdegeneracycount{\exclaw}$ has to be
  at least one, because otherwise $\Cechcube\exclaw$ would be an
  intersection $S$-cube which is not a counterexample by assumption.
  Choose an $s\in S$ such that $\ff_s$ is not injective and write
  $S'\coloneqq S\setminus{\set s}$.
  Since $\ff_s$ is primitive, it is of the form
  \begin{equation*}
    \ff_s= \Id_{\numD{\iin-1}}%
    \join\left(\preimage \ff \iin\to\set{\iin}\right)%
    \join\Id_{\numDminus{n}{\iin}}.%
  \end{equation*}
  for some $\iin\in \n$.
  Denote by $L$, $A$ and $R$ the $S$-claws obtained by restricting
  the $S$-claw $\exclaw$ to $\numD{\iin-1}$, $\set{\iin}$ and
  $\numDminus{n}{\iin}$, respectively.
  Hence we have $\exclaw = L\join A\join R$.
  Denote by $L'$ and $R'$ the $S'$-claws induced from $L$ and $R$,
  respectively.
  Since the restriction of $\ff_s$ to both $\numD{\iin-1}$ and
  $\numDminus{n}{\iin}$ is the identity,
  the edges
  \begin{equation*}
    \Cechcubeindirection Ls\colon\simplex 1\lra\Fun(\powersetop{S'},\Delta)
    \intxt{and}
    \Cechcubeindirection Rs\colon\simplex 1\lra\Fun(\powersetop{S'},\Delta),
  \end{equation*}
  corresponding to the \Cech{} cubes
  $\Cechcube L$ and $\Cechcube R$, are the identity
  on the objects $\Cechcube{L'}$ and $\Cechcube{R'}$ of
  $\Fun(\powersetop{S'},\Delta)$, respectively.
  Denote by
  $\const\colon\Delta\to\Fun(\powersetop{S'},\Delta)$
  the constant-diagram functor and define a cosimplicial object $Y$
  in $\Fun(\powersetop{S'},\Delta)$ by
  \begin{equation*}
    Y\colon
    \Delta
    \xra\const
    \Fun(\powersetop{S'},\Delta)
    \xra{\Cechcube{L'}\join(\blank)\join\Cechcube{R'}}
    \Fun(\powersetop{S'},\Delta)
  \end{equation*}
  Denote by $\calY$ the simplicial object
  \begin{equation*}
    \calY\colon
    \Dop \xra{Y^\op}\Fun(\powersetop{S'},\Delta)^\op
    =\Fun(\powerset{S'},\Dop)
    \xra{\localX\circ\blank}\Fun(\powerset{S'},\localC)
  \end{equation*}
  and by
  \begin{equation*}
    p\colon\Fun(\powerset {S'}, \localC)\lra \Fun(\powersetne {S'},\localC)
  \end{equation*}
  the Cartesian fibration of
  \autoref{lem:cartesiancubesasedges}.
  Observe, that the value of $Y$ at the (active) edge
  $\ff_s\cap{\set \iin}\colon(\preimage{\ff_s}{\iin}\to \set{\iin})$
  is precisely the edge
  $\Cechcubeindirection{\exclaw}s$ in $\Fun(\powerset{S'},\Delta)$
  associated to the \Cech{} cube $\Cechcube\exclaw$.
  By \autoref{lem:cartesiancubesasedges}, the simplicial object $\localX$
  sends the cube $\Cechcube\exclaw$ to a
  Cartesian cube if and only if the edge $\calY(\ff_s\cap{\set\iin})$
  is $p$-Cartesian.

  To complete the proof we set up an application
  of the key lemma (\autoref{lem:activetozero}) to show that this edge
  $\calY(\ff_s\cap{\set\iin})$ is $p$-Cartesian,
  so that the cube $\Cechcube\exclaw$ was not a counterexample after all.
  Let $m\geq 1$ and consider the $S$-claw
  $\exclaw^m=\tupleP{\ff^m_{s'}}{{s'}\in S}$
  on $\numD{\iin-1}\join\m\join{\numDminus{n}{\iin}}$ given by
  \begin{equation*}
    \ff^m_{s'}\coloneqq
    \left(\ff_{s'}\cap {\numD{\iin-1}}\right)
    \join
    \Id_{\m}
    \join
    \left(\ff_{s'}\cap{\numDminus{n}{\iin}}\right)
  \end{equation*}
  for all $s'\neq s$ and by
  \begin{equation*}
    \ff^m_s\coloneqq\Id_{\numD{\iin-1}}%
    \join \left(\canactivemap{m}\colon\numD 1\to\m \right)%
    \join\Id_{\numDminus{n}{\iin}}.%
  \end{equation*}
  It is clear that the $S$-claw $\exclaw^m$
  inherits compatibility from $\exclaw$ and that the \Cech{} cube
  $\Cechcube{\exclaw^m}$ corresponds precisely to the edge
  \begin{equation*}
    Y(\canactivemap{m})\colon
    \simplex 1\xra{\canactivemap m}\Delta\xra{Y}\Fun(\powersetop {S'},\Delta).
  \end{equation*}
  For every $s'\in S\setminus{\set s}$,
  the map $\ff^m_{s'}$ is injective if and only if $\ff_{s'}$ is injective.
  Furthermore, the map $\ff^m_{s}$ is injective
  (this is where we use here we use $m\neq 0$);
  hence the number $\localdegeneracycount {\exclaw^m}$ is smaller than
  $\localdegeneracycount{\exclaw}$.
  By the minimality assumption on the counterexample $\exclaw$,
  we conclude that the simplicial object $\localX$
  sends the \Cech{} cube $\Cechcube{\exclaw^m}$ to a Cartesian cube.
  By \autoref{lem:cartesiancubesasedges} this translates to the fact that the
  corresponding edge
  $\localX\circ\Cechcubeindirection{\exclaw^m}s=\calY(\canactivemap{m})$ in
  $\Fun(\powerset{S'},\localC)$ is $p$-Cartesian.
  Finally, we apply the key lemma
  (\autoref{lem:activetozero}) to the Cartesian fibration $p$
  and the simplicial object $\calY$ to deduce that $\calY$ sends all active maps
  in $\Delta$ to $p$-cartesian edges;
  in particular this is true for the active map
  $\ff_s\cap{\set\iin}\colon\preimage{\ff_s}i\to \set{i}$.
  This completes the proof.
\end{Prf}

\subsection{Proof of the key lemma}
\label{sec:proof-key-lemma}

\begin{Cstr}
  \label{cstr:ractlact}
  Via the functor
  \begin{equation*}
    J\lmapsto J\disjunion\set\infty
  \end{equation*}
  we identify the augmented simplex category $\Deltaug$
  with the wide subcategory
  $\Drstr\subset\Dract$ spanned by the right strict morphisms.
  For every right active morphism $f\colon\m\to\n$ in $\Delta$
  we define a left active morphism $\readbr f\colon\n\to\m$ by the formula
  \begin{equation*}
    \readbr f\colon j\lmapsto \min \invimage{f}\set{j,\dots,n}.
  \end{equation*}
  For every left active morphism $g\colon\n\to\m$ in $\Delta$
  we define a left active morphism $\readbl g\colon\m\to\n$ by the formula
  \begin{equation*}
    \readbl g\colon i\lmapsto \max \invimage{g}\set{0,\dots,i}.
    \qedhere
  \end{equation*}
\end{Cstr}

\begin{Lem}[Joyal duality]
  \label{lem:augisact}
  The assignments $f\mapsto \readbr f$ and $g\mapsto \readbl g$ of
  \autoref{cstr:ractlact}
  are mutually inverse and assemble to an isomorphism of categories
  \begin{equation*}
    \Dract\xlra{\cong}\Dlactop
  \end{equation*}
  (given by the identity on objects) which restricts to an isomorphism
  \begin{equation*}
    \Deltaug\cong\Drstr\xlra{\cong}\Dactop.
    \qedhere
  \end{equation*}
\end{Lem}

\begin{Prf}
  This is a straightforward calculation.
\end{Prf}

The category $\Dact$ has an initial object $\numD 1$ and a terminal object $\0$
which, under the identification $\Deltaug\cong\Dactop$
of \autoref{lem:augisact} correspond
to the objects $\0$ and $\emptyset$ of $\Deltaug$, respectively.

\begin{Lem}
  \label{lem:activeislimit}
  Let $\xdc$ be a simplicial object in any $\infty$-category $\localC$.
  Then the restriction of $\localX$ to the subcategory $\Dactop\subset\Dop$
  is a limit cone.
\end{Lem}

\begin{Prf}
  Lemma~6.1.3.16 in \cite{Lurie2009} states (after passing to opposite
  categories) that every augmented cosimplicial object
  $\Deltaug\cong\Drstr\to\localC$ which extends to a diagram $\Dract\to\localC$
  is automatically a limit diagram. Hence by~\autoref{lem:augisact} every
  diagram $\Dlactop\to\localC$ and, \latin{a fortiori}, every simplicial object
  $\Dop\to\localC$ restricts to a limit diagram $\Dactop\to\localC$.
\end{Prf}

\begin{Prf}[of they key lemma, \autoref{lem:activetozero}]
  \newcommand{\Xact}{\localX^{\mathrm{act}}}
  Denote by $\Xact$ the restriction of $\localX$ to $\Dact$.
  Denote by $\Dactgeq 1$ the full subcategory
  of $\Dact$ spanned by the objects $\m$ with $m\geq 1$.
  Applying \autoref{lem:activeislimit} twice
  we deduce that $\Xact$ and $p\circ\Xact$ are limit cones;
  it follows from \cite[Proposition~4.3.1.5]{Lurie2009}
  that $\Xact$ is also a $p$-limit cone,
  \ie a right $p$-Kan extensions of its restriction to $\Dactgeqop 1$.
  Since the object $\numD 1\in\Dact$ is initial,
  the assumption of \autoref{lem:activetozero} expresses precisely that
  the restriction of $\Xact$ to $\Dactgeqop 1$
  is the right $p$-Kan extension
  of its restriction to $\set{\numD 1}\subset\Dact$.
  We conclude by transitivity of $p$-Kan extensions
  \cite[Proposition~4.3.2.8]{Lurie2009}
  that $\Xact$ is a right $p$-Kan extension
  of its restriction to $\set{\numD 1}$,
  which implies by the pointwise formula at $\0\in\Dact$
  that the edge $\localX(\canactivemap{0}\colon\numD 1\to\0)$ is $p$-Cartesian.
  For every active map $\alpha\colon\m\to\n$ in $\Delta$ we have
  $\canactivemap{n}\circ\alpha=\canactivemap{m}$
  and we already know that the edges
  $\localX(\canactivemap{n})$ and $\localX(\canactivemap{m})$
  are $p$-Cartesian;
  it follows by the left cancellation property of $p$-Cartesian edges
  \cite[Proposition~2.4.1.7]{Lurie2009}
  that the edge $\localX(\alpha)$ is also $p$-Cartesian.
\end{Prf}

\section{Higher Segal conditions}
\label{sec:higher-Segal}

In this last section, we explain the relationship between
the higher Segal spaces of \nameDK{}
and ($\Delta$- and $\Lambda$-)excisive simplicial objects.

\subsection{Higher Segal covers}
\label{sec:higher-segal-covers}

Fix a positive natural number $\kk\geq 1$.
Given a subset $\II\subseteq \n$,
a \introduce{gap} of $\II$ (with $\n$ implicit)
is an element $x\in \n$ with $x\notin\II$.
A gap $x$ of $\II\subseteq\n$ is called
\introduce{even}
if the cardinaity
$\card{\{y\in \II\mid x<y\}}$
is even.
A subset $\II\subseteq \n$ is called \introduce{even}
if all its gaps are even.
Note that even subsets $\II\subseteq \n$ of cardinality $2\kk$
are precisely those which
can be written as a disjoint union of the form
\begin{equation*}
  \II=\bigdisjunion_{i=1}^\kk\{x_i,x_i+1\},
\end{equation*}
with $0\leq x_1<x_1+1<x_2<\dots<x_{k-1}+1<x_k<x_{k}+1\leq n$.

\begin{Def}
  For each $n\geq 2\kk$, the \introduce{\loweroddSegal{\kk} cover}
  on $\n\in \Delta$ is defined as follows:
  \begin{equation*}
    \lSegalclaw\kk n\coloneqq
    \setP{\II\subset\n}
    {\II \text{ even with of cardinality} \card{\II}=2\kk}\clawon\n
    \qedhere
  \end{equation*}
\end{Def}

Observe that the \loweroddSegal{\kk} covers
are precisely the canonical \roughly{good $\kk$-covers}
described in \autoref{sec:manifold-calc-Delta}.
The first \loweroddSegal{\kk} cover
$\lSegalclaw{\kk}{n}$,
\ie the one for $n=2\kk$,
is the unique nondegenerate compatible $\numD\kk$-cover on $\n$.
As $n$ grows bigger,
the behavior of \loweroddSegal{\kk} covers on $\n$
and nondegenerate compatible $\numD\kk$-covers on $\n$
diverges dramatically:
In the first case the number of prongs increasingly rapidly with $\n$,
but each subset of $\n$ remains of constant size $2k$;
in the second case it is the number of prongs ($\kk+1$) that stays constant,
while most of the subsets appearing
in a compatible $\numD\kk$-cover are large.
This dichotomy should remind the reader of the analogous
behavior of $\coveragegood{\kk}$ and $\coveragepoly{\kk}$
described in \autoref{sec:manifold-calc-Delta}:
\begin{itemize}
\item
  Good $\kk$-covers of a manifold typically consist of
  a large number of open subsets;
  however, each of these subsets is simple and small
  (just a disjoint union of at most $\kk$ balls)
\item
  The open covers in $\coveragepoly{\kk}$
  always contain exactly $\kk+1$ open subsets $M\setminus A_i$;
  however, each of these open subsets is usually big and complicated.
\end{itemize}

\begin{Expl}
  The following is a depiction of the first two
  lower $3$-Segal covers:
  \begin{equation}
    \begin{array}{*{5}c}
      0 & 1 & 2 & 3 & 4\\
      \hline
     \xxO&\xxX&\xxX&\xxX&\xxX\\
     \xxX&\xxX&\xxO&\xxX&\xxX\\
     \xxX&\xxX&\xxX&\xxX&\xxO
    \end{array}
    \intxt{and}
    \begin{array}{*{6}c}
      0 & 1 & 2 & 3 & 4 & 5\\
      \hline
     \xxO&\xxO&\xxX&\xxX&\xxX&\xxX\\
     \xxO&\xxX&\xxX&\xxO&\xxX&\xxX\\
     \xxX&\xxX&\xxO&\xxO&\xxX&\xxX\\
     \xxO&\xxX&\xxX&\xxX&\xxX&\xxO\\
     \xxX&\xxX&\xxO&\xxX&\xxX&\xxO\\
     \xxX&\xxX&\xxX&\xxX&\xxO&\xxO
    \end{array}
  \end{equation}
  Observe that the left cover is the unique nondegenerate
  compatible $\numD 2$-cover on $\numD 4=\numD{2\kk}$.
\end{Expl}

We now come to the definition of higher Segal objects.
The definition we will use is not the original one,
but rather a reformulation called the
\buzzword{path space criterion}~\cite[Proposition 2.7]{Poguntke2017}.
\begin{Def}
  A simplicial object $\xdc$ is called
  \begin{itemize}
  \item
    \introduce{\loweroddSegal{\kk}} if,
    for each $n\geq 2\kk$,
    it satisfies descent with respect to the \loweroddSegal{\kk} cover
    $\lSegalclaw\kk n$;
  \item
    \introduce{lower \evenSegal{\kk}} if
    the left path object $\lpath\localX$ is \loweroddSegal{\kk};
  \item
    \introduce{upper \evenSegal{\kk}} if
    the right path object $\rpath\localX$ is \loweroddSegal{\kk};
  \item
    \introduce{\evenSegal{\kk}} if
    $\localX$ is both lower and upper \evenSegal{\kk}.
    \qedhere
  \end{itemize}
\end{Def}

\subsection{Segal = polynomial = weakly excisive}

We come now to the main result of this article,
the comparison of higher Segal conditions and weak excision.
The key ingredient is the following theorem,
which identifies the hierarchy of lower odd Segal objects
with the hierarchy of polynomial functors.

\begin{Thm}
  \label{thm:Segal-poly}
  Let $\localC$ be an \inftycat{} with finite limits.
  The \loweroddSegal{\kk} objects in $\localC$
  are precisely the polynomial functors $\Dop\to\localC$
  of degree $\leq \kk$.
\end{Thm}

Before we prove \autoref{thm:Segal-poly},
we use it to deduce our main theorem.

\begin{Thm}
  \label{thm:main}
  A simplicial object in an $\infty$-category with finite limits is
  \begin{enumerate}
  \item
    \label{item:Segal-exc}
    \loweroddSegal{\kk}
    if and only if
    it is weakly $\kk$-excisive.
  \item
    lower \evenSegal{\kk}
    if and only if
    it is lower weakly $\kk$-$\Lambda$-excisive.
  \item
    upper \evenSegal{\kk}
    if and only if
    it is upper weakly $\kk$-$\Lambda$-excisive.
  \item
    \evenSegal{\kk}
    if and only if
    it is weakly $\kk$-$\Lambda$-excisive.
    \qedhere
  \end{enumerate}
\end{Thm}

 \begin{Prf}[of \autoref{thm:main}]
   In \autoref{thm:unitality-of-cubes}
   we have seen that a functor $\Dop\to \localC$
   is polynomial of degree $\leq\kk$
   if and only if it is weakly $\kk$-excisive;
   thus part \autoref{item:Segal-exc}
   is an immediate consequence of \autoref{thm:Segal-poly}.
   The rest of \autoref{thm:main} then follows immediately
   from the path space criterion for weak $\Lambda$-excision
   (\autoref{cor:path-space-criterion-cubes}).
\end{Prf}

Recall that a cyclic object $\Lambda^\op\to\localC$
is defined to be \evenSegal{\kk}
if the underlying simplicial object
$\Dop\to\Lambda^\op\to\localC$
is \evenSegal{\kk}.

\begin{Cor}
  \label{cor:cyclicmain}
  A cyclic object in an $\infty$-category with finite limits
  is \evenSegal{\kk} if and only if it is weakly $\kk$-excisive.
\end{Cor}

\begin{Prf}
  \autoref{cor:cyclicmain} follows directly from
  \autoref{thm:main} and \autoref{rem:cyclic-excisive}.
\end{Prf}

We now give the proof of \autoref{thm:Segal-poly}.

\begin{Prf}[of \autoref{thm:Segal-poly}]
  \newcommand\localspF[2]{\exclaw_{#2}^{#1}}
  \newcommand{\inprfm}{m}
  \newcommand\localspf[2]{\II_{#2}^{#1}}
  Fix a simplicial object $\xdc$
  in an $\infty$-category $\localC$ with finite limits.
  By the characterization of
  strongly biCartesian intersection cubes in $\Delta$
  (\autoref{cor:biCart=comp})
  we only need to show that
  $\localX$ satisfies descent
  with respect to all \loweroddSegal{\kk} covers
  if and only if
  $\localX$ satisfies descent with respect to all
  compatible $\kk$-covers.
  In view of \autoref{lem:onebiCartforall},
  we only have to relate, for each $n\geq2\kk$,
  the \loweroddSegal{\kk} cover to \emph{one}
  nondegenerate compatible $\kk$-cover.
  For each $n\geq 2\kk$ and each $j\in\set{-1,0,\dots,\kk}$,
  we define a cover $\localspF nj\clawon\n$
  (with the $\kk$ left implicit since it is fixed throughout the proof)
  to consist of the following subsets of $\n$:

  \begin{itemize}
  \item
    $\localspf ni\coloneqq \n\setminus \set{2i}$ for $i=0,\dots,j$
  \item
    those $\II\in \lSegalclaw n\kk$ that satisfy
    $\numD{2j}=\set{0,1,\dots,2j}\subset \II$.
  \end{itemize}


  Clearly $\localspF n{-1}$ is nothing but
  the \loweroddSegal{\kk} cover $\lSegalclaw n\kk\clawon\n$.
  Moreover, we have a chain of refinements
  \begin{equation}
    \label{eq:interpolationrefinement}
    \lSegalclaw n\kk
    =\localspF n{-1}
    \refines\localspF n0
    \refines\dots
    \refines\localspF n\kk
  \end{equation}
  because every $\II\in\lSegalclaw n\kk$
  with $\numD{2(j-1)}\subset \II$
  must either satisfy $\numD{2j}\subset \II$ or $2j\notin \II$.
  The last cover
  $\localspF n\kk=\tupleP{\localspf n i}{i\in\kk}$
  in the refinement~\eqref{eq:interpolationrefinement}
  is a nondegenerate compatible $\numD\kk$-claw;
  in this sense, the chain~\eqref{eq:interpolationrefinement}
  is an interpolation between the Segal condition and
  the descent condition with respect to the family
  $\setP{\localspF n \kk}{n\geq 2\kk}$
  of nondegenerate compatible $\numD\kk$-covers in $\Delta$.

  We establish the following two facts:
  \begin{enumerate}
  \item
    \label{item:basecasec}
    If $n=2\kk$ then the chain~\eqref{eq:interpolationrefinement} of
    refinements collapses, \ie we have
    \begin{equation*}
      \lSegalclaw {2\kk}\kk
      =\localspF {2\kk}{-1}
      =\localspF {2\kk}0
      =\dots
      =\localspF {2\kk}\kk.
    \end{equation*}
  \item
    \label{item:indcasew}
    For every $n>2\kk$ and every $j=0,\dots,n$
    the refinement
    $\localspF n{j-1} \refines\localspF nj$
    is $\localX$-local provided that the cover
    $\localspF {n-1}{j-1}\clawon\numD{n-1}$
    is $\localX$-local.
  \end{enumerate}
  Fact~\ref{item:basecasec} is immediate from the definition.
  For each $j=0,\dots,\kk$ we have
  $\localspF n j =\reducedclaw{\localspF n {j-1}\cup\set{\localspf nj}}$
  and the cover
  $\localspF n{j-1}\cap\localspf n{j}\clawon\localspf n{j}$
  is easily seen to be isomorphic
  (under the unique isomorphism $\localspf n j\cong\numD{n-1}$)
  to the cover $\localspF {n-1}{j-1}\clawon \numD{n-1}$;
  hence fact~\ref{item:indcasew} follows from
  \autoref{lem:one-step-refinement}.

  By a straightforward inductive argument,
  facts~\ref{item:basecasec} and~\ref{item:indcasew} imply
  that the following three conditions are equivalent:
  \begin{itemize}
  \item
    For all $n\geq 2\kk$,
    the cover
    $\lSegalclaw{n}{\kk}=\localspF n{-1}\clawon\n$
    is $\localX$-local.
  \item
    For all $n\geq 2\kk$ and all $j=-1,\dots,\kk$,
    the cover
    $\lSegalclaw{n}{\kk}=\localspF n{j}\clawon\n$
    is $\localX$-local.
  \item
    For all $n\geq 2\kk$,
    the (nondegenerate, compatible, $\numD\kk$-pronged) cover
    $\localspF n\kk\clawon\n$
    is $\localX$-local.
  \end{itemize}
  We have therefore related the Segal conditions to
  one hierarchy of descent conditions
  with respect to nondegenerate compatible $\numD\kk$-covers;
  \autoref{lem:onebiCartforall} precisely states that this is enough,
  hence the proof is concluded.
\end{Prf}

\subsection{Triviality bounds for higher Segal objects}

Let $\xdc$ be a lower or upper $\dd$-Segal object in $\localC$.
Since for each $m>\dd$ the $\dd$-Segal conditions express
the value $\localX_m$ as a cubical limit
of the values $\localX_{n}$ with $n\leq \dd$,
it is obvious that $\localX$ is trivial
(\ie $\localX_n$ is a terminal object in $\localC$ for each $\n\in\Delta$)
as soon as $\localX$ is trivial when restricted to $\Deltaleq{\dd}$.
From the comparison with weak excision we can deduce
the following sharper bounds:

\begin{Prop}
  \label{prop:triviality-bounds}
  Fix $\dd\geq 2$ and
  let $\xdc$ be a lower or upper $\dd$-Segal object
  in an \inftycat{} $\localC$ with finite limits.
  If $\localX$ is trivial when restricted to $\Deltall{\dd}$
  then $\localX$ is trivial.
\end{Prop}

\begin{Rem}
  Since not every monoid is trivial,
  it is not true that a lower $1$-Segal object
  (\ie a Segal object in the sense of Rezk) is trivial as soon as
  its restriction to $\Deltaleq{0}$ is trivial.
  Hence the assumption $\dd\geq 2$ in \autoref{prop:triviality-bounds}
  is necessary.
\end{Rem}

\begin{Prf}[of \autoref{prop:triviality-bounds}]
  First, we prove the case of lower odd Segal objects.
  Let $\kk\geq 2$ and assume that $\xdc$ is \loweroddSegal{\kk}
  and trivial on $\Deltall{2\kk-1}$.
  It suffices to show that $\localX_{\numD{2\kk-1}}$ is trivial.
  Consider the following compatible $\numD\kk$-claw $\exclaw$
  on $\numD{2\kk-2}$:
  \begin{equation}
    \label{eq:4}
    \begin{array}{c|*{9}c}
             & 0 & 1  & 2 & 3 & 4 & \cdots & 2\kk-4 & 2\kk-3 & 2\kk-2\\
      \hline
      0      &\xxX&\xx2&\xxX&\xxX&\xxX& \cdots &\xxX&\xxX&\xxX\\
      1      &\xxO&\xxX&\xxX&\xxX&\xxX& \cdots &\xxX&\xxX&\xxX\\
      2      &\xxX&\xxX&\xxO&\xxX&\xxX& \cdots &\xxX&\xxX&\xxX\\
      3      &\xxX&\xxX&\xxX&\xxX&\xxO& \cdots &\xxX&\xxX&\xxX\\
      \vdots\\
      \kk-1  &\xxX&\xxX&\xxX&\xxX&\xxX& \cdots &\xxO&\xxX&\xxX\\
      \kk    &\xxX&\xxX&\xxX&\xxX&\xxX& \cdots &\xxX&\xxX&\xxO\\
    \end{array}
  \end{equation}
  The corresponding biCartesian \Cech{} cube
  $\Cechcube\exclaw\colon\powersetop{\numD\kk}\to\Delta$
  satisfies $\Cechcube\exclaw(\set{0})\cong\numD{2\kk-1}$
  and $\Cechcube\exclaw(T)\in \Deltall{2\kk-1}$ for all $T\neq\set{0}$.
  It follows that the $\numD\kk$-cube
  $\localX\circ\Cechcube\exclaw$
  sends every $T\subseteq\numD\kk$, except possibly $T=\set{0},$
  to a terminal object in $\localC$.
  Since $\localX$ is weakly $\kk$-excisive by \autoref{thm:main},
  this cube in $\localC$ is Cartesian.
  It then follows that we have a Cartesian square
  \begin{equation*}
    \cdsquareNA[C]
    {(\localX\circ\Cechcube\exclaw)(\emptyset)}
    {\lim\limits_{\substack{
          0\notin T\subseteq \numD\kk\\
          \emptyset\neq T
        }}
      (\localX\circ\Cechcube\exclaw)(T)}
    {(\localX\circ\Cechcube\exclaw)(\set{0})}
    {\lim\limits_{\substack{
          0\in T\subseteq \numD\kk\\
          \set{0}\neq T
        }}
      (\localX\circ\Cechcube\exclaw)(T)}
  \end{equation*}
  in $\localC$,
  where all but the lower left corner are trivial;
  we conclude that
  $(\localX\circ\Cechcube\exclaw)(\set{0})\simeq\localX_{\numD{2\kk-1}}$
  is also trivial.

  If $\dd=2\kk$ is even with $\kk\geq 1$ then the same proof works for
  lower or upper \evenSegal{\kk} objects
  by considering instead of \eqref{eq:4}
  the left active compatible $\kk$-claw
  \begin{equation*}
    \begin{array}{c|*{8}c}
             & 0 & 1  & 2 & 3 & \cdots & 2\kk-3 & 2\kk-2 & 2\kk-1\\
      \hline
      0      &\xx2&\xxX&\xxX&\xxX& \cdots &\xxX&\xxX&\xxX\\
      1      &\xxX&\xxO&\xxX&\xxX& \cdots &\xxX&\xxX&\xxX\\
      2      &\xxX&\xxX&\xxX&\xxO& \cdots &\xxX&\xxX&\xxX\\
      \vdots\\
      \kk-1  &\xxX&\xxX&\xxX&\xxX& \cdots &\xxO&\xxX&\xxX\\
      \kk    &\xxX&\xxX&\xxX&\xxX& \cdots &\xxX&\xxX&\xxO\\
    \end{array}
  \end{equation*}
  on $\numD{2\kk-1}$ or its obvious right active analog.

  Recall from~\cite[Proposition~2.7]{Poguntke2017}
  that a simplicial object is upper ($2\kk+1$)-Segal
  if and only if its left path object is upper \evenSegal{\kk}
  (or, equivalently, if its right path object is lower \evenSegal{\kk});
  the result for upper odd Segal objects thus follows
  immediately from the one for (lower or upper) even Segal objects.
\end{Prf}

It is not known to \theauthor{} if the bounds in
\autoref{prop:triviality-bounds} are sharp.
More precisely,
\theauthor{} does not know the answer to the following question,
which remains to be investigated in future work:

\begin{Qstn}
  Let $\kk\geq 1$ and
  let $\localX$ be a simplicial object which is
  \loweroddSegal{\kk},
  or upper \evenSegal{\kk}
  or lower \evenSegal{\kk}.
  If $\localX$ is trivial when restricted to $\Deltaleq\kk$,
  does it follow that $\localX$ is trivial?
\end{Qstn}

\newpage
\bibliographystyle{amsalpha}
\bibliography{library}
\end{document}